\documentclass{jnmp}
\usepackage{amsmath,amssymb}
\usepackage[latin1]{inputenc}  
\JNMPnumberwithin{equation}{section}
\resetfootnoterule

\newtheorem{theorem}{Theorem}[section]
\newtheorem{lemma}{Lemma}[section]
\newtheorem{corollary}{Corollary}[section]
\newtheorem{proposition}{Proposition}[section]

\theoremstyle{definition}
\newtheorem{definition}{Definition}[section]
\newtheorem{remark}{Remark}[section]
\newtheorem*{notation}{Notation}

\begin{document}

\renewcommand{\evenhead}{M Sep\'{u}lveda and O Vera}
\renewcommand{\oddhead}{Analycity for the coupled system of KdV equations}

\thispagestyle{empty}

\FirstPageHead{*}{*}{20**}{\pageref{firstpage}--\pageref{lastpage}}{Article}

\copyrightnote{200*}{M Sep\'{u}lveda and O Vera}

\Name{Analycity and smoothing effect for the coupled system of
equations of Korteweg - de Vries type with a single point
singularity}

\label{firstpage}

\Author{Mauricio Sep\'{u}lveda~$^a$ and Octavio Vera
Villagr\a'{a}n~$^b$}

\Address{$^a$ Departamento de Ingenier\'{\i}a Matem\a'{a}tica, Universidad de
Concepci\a'{o}n, Casilla 160-C, Concepci\a'{o}n, Chile. \\
~~E-mail: mauricio@ing-mat.udec.cl\\[10pt]
$^b$ Departamento de Matem\a'{a}tica, Universidad
del B\a'{\i}o-B\a'{\i}o, Collao 1202, Casilla 5-C, Concepci\a'{o}n,
Chile. \\
~~E-mail: overa@ubiobio.cl}

\Date{Received Month *, 200*; Accepted in Revised Form Month *, 200*}

\begin{abstract}
\noindent We study that a solution of the initial value problem
associated for the coupled system of equations of Korteweg - de
Vries type which appears as a model to describe the strong
interaction of weakly nonlinear long waves, has  analyticity in
time and smoothing effect up to real analyticity if the initial data
only has a single point singularity at $x=0.$
\end{abstract}

\renewcommand{\theequation}{\thesection.\arabic{equation}}
\setcounter{equation}{0}\section{Introduction} We consider the
following coupled system of equations of Korteweg - de Vries type
\begin{eqnarray}
\label{e101}&  &\widetilde{u}_{t} + \widetilde{u}_{xxx} +
a_{3}\,\widetilde{v}_{xxx} + \widetilde{u}\,\widetilde{u}_{x} +
a_{1}\,\widetilde{v}\,\widetilde{v}_{x}
 + a_{2}\,(\widetilde{u}\,\widetilde{v})_{x} = 0,\quad
 x,\,t\in\mathbb{R} \\
\label{e102}&  & b_{1}\,\widetilde{v}_{t} + \widetilde{v}_{xxx} +
b_{2}\,a_{3}\,\widetilde{u}_{xxx} + \widetilde{v}\,\widetilde{v}_{x}
+ b_{2}\,a_{2}\,\widetilde{u}\,\widetilde{u}_{x}
+ b_{2}\,a_{1}\,(\widetilde{u}\,\widetilde{v})_{x} = 0,  \\
\label{e103}&  &\widetilde{u}(x,\,0)=\widetilde{u}_{0}(x),\quad
\widetilde{v}(x,\,0)=\widetilde{v}_{0}(x).
\end{eqnarray}
where
$\widetilde{u}=\widetilde{u}(x,\,t),\,\widetilde{v}=\widetilde{v}(x,\,t)$
are real-valued functions of the variables $x$ and $t$ and
$a_{1},\,a_{2},\,a_{3},\,b_{1},\,b_{2}$ are real constants with
$b_{1}>0$ and $b_{2}>0.$  The original coupled system is
\begin{eqnarray}
\label{e104}&  & \widetilde{u}_{t} + \widetilde{u}_{xxx} +
a_{3}\,\widetilde{v}_{xxx} + \widetilde{u}^{p}\,\widetilde{u}_{x} +
a_{1}\,\widetilde{v}^{p}\,\widetilde{v}_{x}
 + a_{2}\,(\widetilde{u}^{p}\,\widetilde{v})_{x} = 0,\quad
 x,\,t\in\mathbb{R}   \\
\label{e105}&  & b_{1}\,\widetilde{v}_{t} + \widetilde{v}_{xxx} +
b_{2}\,a_{3}\,\widetilde{u}_{xxx} +
\widetilde{v}^{p}\,\widetilde{v}_{x} +
b_{2}\,a_{2}\,\widetilde{u}^{p}\,\widetilde{u}_{x} + b_{2}\,a_{1}\,
(\widetilde{u}\,\widetilde{v}^{p})_{x}  =  0 \\
\label{e106}&  & \widetilde{u}(x,\,0)=\widetilde{u}_{0}(x),\quad
\widetilde{v}(x,\,0)=\widetilde{v}_{0}(x)
\end{eqnarray}
where
$\widetilde{u}=\widetilde{u}(x,\,t),\,\widetilde{v}=\widetilde{v}(x,\,t)$
are real-valued functions of the variables $x$ and $t$ and
$a_{1},\,a_{2},\,a_{3},\,b_{1},\,b_{2}$ are real constants with
$b_{1}>0$ and $b_{2}>0.$ The power $p$ is an integer larger than or
equal to one. The system \eqref{e104}-\eqref{e106} has the structure
of a pair of Korteweg - de Vries equations coupled through both
dispersive and nonlinear effects. In the case $p=1,$ the system
\eqref{e104}-\eqref{e106} was derived by Gear and Grimshaw \cite{ge}
as a model to describe the strong interaction of weakly nonlinear,
long waves. Mathematical results on the system
\eqref{e104}-\eqref{e106} were given by J. Bona {\it et al.}
\cite{bo3}. They proved that \eqref{e104}-\eqref{e106} is globally
well posed in $H^{s}(\mathbb{R})\times H^{s}(\mathbb{R})$ for any
$s\geq 1$ provided $|a_{3}|<1/\sqrt {b_{2}}.$ The system
\eqref{e104}-\eqref{e106} has been intensively studied by several
authors (see \cite{bi1,bi2,bo3,da1,thesi1} and the references therein). We have
the following conservation laws
\begin{eqnarray}
\mathbb{E}_{1}(\widetilde{u})=\int_{\mathbb{R}}\widetilde{u}\,dx\quad
,\quad \mathbb{E}_{2}(\widetilde{v})=
\int_{\mathbb{R}}\widetilde{v}\,dx\quad ,\quad \mathbb{E}
_{3}(\widetilde{u},\,\widetilde{v})=\int_{\mathbb{R}}(b_{2}\widetilde{u}^{2}
+ b_{1}\widetilde{v}^{2})dx
\end{eqnarray}
The time-invariance of the functionals $\mathbb{E}_{1}$ and
$\mathbb{E}_{2}$ expresses the property that the mass of each mode
separately is conserved during interaction, while that of
$\mathbb{E}_{3}$ is an expression of the conservation of energy for
the system of two models taken as a whole. The solutions of
\eqref{e104}-\eqref{e106} satisfy an additional conservation law
which is revealed by the time-invariance of the functional
\begin{eqnarray*}
\label{e108}\mathbb{E}_{4}=\int_{\mathbb{R}}\left(b_{2}\,\widetilde{u}_{x}^{2}
+ \widetilde{v}_{x}^{2} +
2b_{2}a_{3}\widetilde{u}_{x}\widetilde{v}_{x} - b_{2}\frac
{\widetilde{u}^{3}}{3} - b_{2}a_{2}\widetilde{u}^{2}\widetilde{v} -
b_{2}a_{2}\widetilde{u}^{2}\widetilde{v} -
b_{2}a_{1}\widetilde{u}\widetilde{v}^{2} - \frac
{\widetilde{v}^{3}}{3}\right)dx
\end{eqnarray*}
The functional $\mathbb{E}_{4}$ is a Hamiltonian for the system
\eqref{e104}-\eqref{e106} and if $b_{2}a_{3}^{2}<1,\,\phi _{4}$ will
be seen to provide an a priori estimate for the solutions
$(\widetilde{u},\,\widetilde{v})$ of \eqref{e104}-\eqref{e106} in
the space $H^{1}(\mathbb{R})\times H^{1}(\mathbb{R}).$ Furthermore,
the linearization of \eqref{e101}-\eqref{e103} about the rest state
can be reduced to two, linear Korteweg - de Vries equations by a
process of diagonalization. Using this remark and the smoothing
properties (in both the temporal and spatial variables) for the
linear Korteweg - de Vries derived by Kato \cite{ka4,ka6}, Kenig,
Ponce and Vega \cite{ke4,ke5} it will be shown that
\eqref{e104}-\eqref{e106} is locally well-posed in
$H^{s}(\mathbb{R})\times H^{s}(\mathbb{R})$ for any $s\geq 1$
whenever $\sqrt {b_{2}}a_{3}\neq 1.$ This result was improved by J.
M. Ash {\it et al.} \cite{a0} showing that the system
\eqref{e101}-\eqref{e103} is globally well-posed in
$L^{2}(\mathbb{R})\times L^{2}(\mathbb{R})$ provided that $\sqrt
{b_{2}}a_{3}\neq 1.$ In 2004, F. Linares and M. Panthee \cite{li1}
improve this result showing that the system \eqref{e101}-\eqref{e103}
is locally well-posed in $H^{s}(\mathbb{R})\times H^{s}(\mathbb{R})$
for $s>-3/4$ and globally well-posed in $H^{s}(\mathbb{R})\times
H^{s}(\mathbb{R})$ for $s>-3/10$ under some conditions on the
coefficients, indeed for $a_{3}=0$ and $b_{1}=b_{2}.$  
Following the idea
W. Craig {\it et al.} \cite{cr3}, it is  shown in \cite{thesi1}
that $C^{\infty}$ solutions
$(\widetilde{u}(\,\cdot\,,\,t),\,\widetilde{v}(\,\cdot\,,\,t))$ to
\eqref{e101}-\eqref{e103} are obtained for $t>0$ if the initial data
$(\widetilde{u}(x,\,0),\,\widetilde{v}(x,\,0))$ belong to a suitable
Sobolev space satisfying resonable conditions as
$|x|\rightarrow\infty.$ Since \eqref{e101}-\eqref{e103} is a coupled
system of Korteweg-de Vries equations, it is natural to ask whether it
has a smoothing effect up to real analyticity if the initial data
only has a single point singularity at $x=0$ as the known results for
the scalar case of a single Korteweg -de Vries equation. 
Using the scaling argument we can have
an insight to this question. In this paper our purpose is to prove
the analyticity in time of solutions to \eqref{e101}-\eqref{e103}
without regularity assumption on the initial data improving those
obtained in \cite{thesi1}. Our main tool is the
generator of dilation $P=3\,t\,\partial_{t} + x\,\partial_{x}.$
which almost commutes with the linear Korteweg-de Vries operator
$L=\partial_{t} +
\partial_{x}.$ Indeed $[L,\,P]=3\,L.$ A typical example of initial
data satisfying the assumption of the above theorem is the Dirac
delta measure, since
$(x^{k}\,\partial_{x})^{k}\delta(x)=(-1)^{k}\,{k!}\,\delta(x).$ The
other example of the data is p. v. $\frac{1}{x},$ where p. v.
denotes the Cauchy principal value. Linear combination of those
distributions with analytic $H^{s}$ data satisfying the assumption
is also possible. In this sense, the Dirac delta measure adding the
soliton initial data can be taken as an initial datum. Using the
operator $K=x\cdot\bigtriangledown + 2\,i\,t\,\partial_{t}$ it was
proved the Gevrey smoothing effect in space variable \cite{de1}.
Indeed, it was shown that, if  the initial data belongs to a Gevrey class
of order 2, then solutions of some nonlinear Schr\"{o}dinger
equations become analytic in the space variable for $t\neq 0.$ For
the Korteweg-de Vries equations version of the generator of dilation
is also useful to study the analyticity in time and the Gevrey effect in
the space variables for solutions \cite{de1}.\\
This paper is organized as follows: In section 2 we have the
reduction of the problem and we outline briefly the notation,
terminology to be used subsequently and results that will be used
several times. In section 3 we prove a theorem of existence and
well-posedness of the solutions. In section 4 we prove the following
theorem:
\begin{theorem} {\it Suppose that the initial data
$(\widetilde{u}_{0},\,\widetilde{v}_{0})\in H^{s}(\mathbb{R})\times
H^{s}(\mathbb{R}),$ $s>-3/4$ and $A_{0},\,A_{1}>0$ such that
\begin{eqnarray}
\label{e109}\sum_{k=0}^{\infty}\frac{A_{0}^{k}}{k!}\,||(x\,\partial_{x})^{k}
\widetilde{u}_{0}||_{H^{s}(\mathbb{R})}<+\infty\qquad :\qquad
\sum_{k=0}^{\infty}\frac{A_{1}^{k}}{k!}\,||(x\,\partial_{x})^{k}
\widetilde{u}_{0}||_{H^{s}(\mathbb{R})}<+\infty.
\end{eqnarray}
Then for some $b\in (1/2,\,7/12),$ there exist
$T=T(||\widetilde{u}_{0}||_{H^{s}(\mathbb{R})},\,||
\widetilde{v}_{0}||_{H^{s}(\mathbb{R})})$ and a unique solution of
\eqref{e101}-\eqref{e103} in a certain time $(-T,\,T)$ and the
solution $(\widetilde{u},\,\widetilde{v})$ is time locally
well-posed, i. e., the solution continuously depends on the initial
data. Moreover, the solution $(\widetilde{u},\,\widetilde{v})$ is
analytic at any point $(x,\,t)\in\mathbb{R}\times\{(-T,\,0)\cup
(0,\,T)\}.$}
\end{theorem}
\begin{corollary} {\it Let $s>-3/4,$ $b\in (1/2,\,7/12).$ Suppose
that the initial data $(\widetilde{u}_{0},\,\widetilde{v}_{0})\in
H^{s}(\mathbb{R})\times H^{s}(\mathbb{R}),$ and $A_{0},\,A_{1}>0$
such that}
\begin{eqnarray}
\label{e110}\sum_{k=0}^{\infty}\frac{A_{0}^{k}}{(k!)^{3}}\,||(x\,\partial_{x})^{k}
\widetilde{u}_{0}||_{H^{s}(\mathbb{R})}<+\infty\qquad :\qquad
\sum_{k=0}^{\infty}\frac{A_{1}^{k}}{(k!)^{3}}\,||(x\,\partial_{x})^{k}
\widetilde{u}_{0}||_{H^{s}(\mathbb{R})}<+\infty.
\end{eqnarray}
{\it Then there exists a unique solution
$(\widetilde{u},\,\widetilde{v})\in
C((-T,\,T),\,H^{s}(\mathbb{R}))\cap X_{b}^{s}\times
C((-T,\,T),\,H^{s}(\mathbb{R}))\cap X_{b}^{s}$ to the coupled system
of Korteweg- de Vries equation \eqref{e101}-\eqref{e103}
 for a certain $(-T,\,T)$ and for any $t\in (-T,\,0)\cup (0,\,T),$
 the pair $(\widetilde{u},\,\widetilde{v})$
 are analytic functions in the space variable and for
 $x\in\mathbb{R},$ $\widetilde{u}(x,\,\cdot\,)$ and
 $\widetilde{v}(x,\,\cdot\,)$ are Gevrey 3 as function of the time
 variable.}
\end{corollary}
\begin{remark} In Theorem 1.1 and Corollary 1.2, the assumption on
 the initial data implies analyticity and Gevrey 3 regularity except
 at the origin respectively. In this sense, those results 
 state that the singularity at the origin immediately disappears
 after $t>0$ or $t<0,$ up to analyticity.
\end{remark}
\begin{remark}
  The crucial part for obtaining a full regularity is
 to gain the $L^{2}(\mathbb{R}^{2})$ regularity of the solutions
 $(u_{k},\,v_{k})$ from the negative order Sobolev space. This part
 is obtained in Proposition 4.1 in Section 4. We utilize a three steps
 recurrence argument for treating the nonlinearity appearing in the
 right hand side of
\begin{eqnarray}
\label{e111}&  & t\,\partial_{x}^{3}u_{k}  =  -\frac{1}{3}\,Pu_{k} +
\frac{1}{3}\,x\,\partial_{x}u_{k} + t\,B_{k}^{1}(u,\,u) +
t\,B_{k}^{2}(v,\,v) + t\,B_{k}^{3}(u,\,v) \\
\label{e112}&  & t\,\partial_{x}^{3}v_{k}  =  -\frac{1}{3}\,Pv_{k} +
\frac{1}{3}\,x\,\partial_{x}v_{k} + t\,C_{k}^{1}(u,\,u) +
t\,C_{k}^{2}(v,\,v) + t\,C_{k}^{3}(u,\,v).
\end{eqnarray}
Then step by step, we obtain the pointwise analytic estimates
\begin{eqnarray}
\label{e113}\sup_{t\in[t_{0} - \epsilon,\,t_{0} +
\epsilon]}||\partial_{t}^{m}\partial_{x}^{l}u||_{H^{1}(x_{0} -
\epsilon,\,x_{0} + \epsilon)}\leq c\,A_{1}^{m + l}\,(m + l)!,\qquad
l,\,m=0,\,1,\,2,\,\ldots \\
\label{e114}\sup_{t\in[t_{0} - \epsilon,\,t_{0} +
\epsilon]}||\partial_{t}^{m}\partial_{x}^{l}v||_{H^{1}(x_{0} -
\epsilon,\,x_{0} + \epsilon)}\leq c\,A_{2}^{m + l}\,(m + l)!,\qquad
l,\,m=0,\,1,\,2,\,\ldots
\end{eqnarray}
Since initially we do not know whether the solution belong to even
$L^{2}(\mathbb{R}^{2})$ we should mention that the local
well-posedness is essentially important for our argument and
therefore it merely satisfies the coupled system equations in the
sense of distribution.
\end{remark}
\renewcommand{\theequation}{\thesection.\arabic{equation}}
\setcounter{equation}{0}\section{Reduction of the Problem and
Preliminary Results} As mentioned in the introduction we consider
the following coupled system of equations of Korteweg - de Vries
type \eqref{e101}-\eqref{e103}.
If $a_{3}=0$ there is no coupling in the
dispersive terms. Let us assume that $a_{3}\neq 0.$ We are interested in
decoupling the dispersive terms in the system
\eqref{e101}-\eqref{e103}. For this, let $a_{3}^{2}\,b_{2}\neq 1.$
We consider the associated linear system
\begin{eqnarray}
\label{e204}W_{t} + A\,W_{xxx} =0,\qquad W(x,\,0)=W_{0}(x)
\end{eqnarray}
where
\[W=\left[ \begin{array}{c}
u\\v \end{array} \right]\qquad ,\qquad A=\left[ \begin{array}{cc}
1 & a_{3}\\
\frac{a_{3}\,b_{2}}{b_{1}} & \frac{1}{b_{1}}
\end{array}
\right].
\]
The eigenvalues of $A$ are given by
\begin{eqnarray}
\label{e205}\alpha_{+} = \frac{1}{2}\left(1 + \frac{1}{b_{1}} +
\sqrt{\left(1 - \frac{1}{b_{1}}\right)^{2} +
\frac{4\,b_{2}\,a_{3}^{2}}{b_{1}}}\right)\\
\label{e206}\alpha_{-} = \frac{1}{2}\left(1 + \frac{1}{b_{1}} -
\sqrt{\left(1 - \frac{1}{b_{1}}\right)^{2} +
\frac{4\,b_{2}\,a_{3}^{2}}{b_{1}}}\right)
\end{eqnarray}
which are distinct since $b_{1}>0,$ $b_{2}>0$ and $a_{3}\neq 0.$ Our
assumption $a_{3}^{2}\,b_{2}\neq 1$ guarantees that
$\alpha_{\pm}\neq 0.$ Thus we can write the system
\eqref{e101}-\eqref{e103} in a matrix form as in \cite{li1}. After
we make the change of scale
\begin{eqnarray*}
\widetilde{u}(x,\,t) =
u(\alpha_{+}^{-1/3}\,x,\,t)\quad\mbox{and}\quad \widetilde{v}(x,\,t)
= v(\alpha_{-}^{-1/3}\,x,\,t).
\end{eqnarray*}
Then we obtain the system of equations
\begin{eqnarray}
\label{e207}&  & u_{t} + u_{xxx} + a\,u\,u_{x} + b\,v\,v_{x}
 + c\;(u\,v)_{x} = 0,\quad
 x,\,t\in\mathbb{R} \\
\label{e208}&  & v_{t} + v_{xxx} + \widetilde{a}\,u\,u_{x} +
\widetilde{b}\,v\,v_{x}
+ \widetilde{c}\;(u\,v)_{x} = 0,  \\
\label{e209}&  & u(x,\,0)=u_{0}(x),\quad v(x,\,0)=v_{0}(x)
\end{eqnarray}
where $a,$ $b,$ $c$ and $\widetilde{a},$ $\widetilde{b},$ $\widetilde{c}\;$
are constant.
\begin{remark} Notice that the nonlinear terms involving the
functions $u$ and $v$ are not evaluated at the same point. Therefore
those terms are not local anymore.
\end{remark}
For $s,\,b\in\mathbb{R}$ define the spaces $X_{b}^{s}$ and $X_{b -
1}^{s}$ to be the completion of the Schwartz space ${\cal
S}(\mathbb{R}^{2})$ with respect to the norms
\begin{eqnarray*}
||u||_{X_{b}^{s}} = \left(\int_{\mathbb{R}}\int_{\mathbb{R}}(1 +
|\tau - \xi^{3}|)^{2b}\,(1
+|\xi|)^{2s}\,|\widehat{u}(\xi,\,\tau)|^{2}\,d\xi\,d\tau\right)^{1/2}
\end{eqnarray*}
and
\begin{eqnarray*}
||u||_{X_{b - 1}^{s}} = \left(\int_{\mathbb{R}}\int_{\mathbb{R}}(1 +
|\tau - \xi^{3}|)^{2(b - 1)}\,(1
+|\xi|)^{2s}\,|\widehat{u}(\xi,\,\tau)|^{2}\,d\xi\,d\tau\right)^{1/2}
\end{eqnarray*}
where $X_{b}^{s}=\{u\in{\cal
S}'(\mathbb{R}^{2}):\;||u||_{X_{b}^{s}}<\infty \}.$ Let ${\cal
F}_{x}$ and ${\cal F}_{x,\,t}$ be the Fourier transform in the $x$
and $(x,\,t)$ variables respectively. The Riesz operator $D_{x}$ is
defined by $D_{x}={\cal F}_{\xi}^{-1}\,|\xi|\,{\cal F}_{x}.$ The
fractional derivative is defined by
\begin{eqnarray*}
<D_{x}>^{s} & = & {\cal F}_{\xi}^{-1}\,<\xi>^{s}\,{\cal
F}_{x} = {\cal F}_{\xi}^{-1}\,(1 + |\xi|^{2})^{s/2}\,{\cal F}_{x}\\
<D_{x,\,t}>^{s} & = & {\cal F}_{\xi,\,\tau}^{-1}\,<|\xi| +
|\tau|>^{s}\,{\cal F}_{x,\,t}
\end{eqnarray*}
For $<\,\cdot\,>=(1 + |\,\cdot\,|^{2})^{1/2},$ we have \\
\\
i) $||\,\cdot\,||_{H^{b}(\mathbb{R}:\,H^{r}(\mathbb{R}))}=
||<D_{t}>^{b}<D_{x}>^{r}\cdot\,||_{L_{x,\,t}^{2}(\mathbb{R}^{2})}.$
\\
ii) $H^{s}(\mathbb{R})=\{u\in {\cal
S}'(\mathbb{R}):\;<D_{x}>^{s}u\in L^{2}(\mathbb{R})\}.$ \\
iii) $||\,\cdot\,||_{H^{s}(\mathbb{R})}=
||<D_{x}>^{s}\cdot\,||_{L^{2}(\mathbb{R})}.$ 
\begin{remark} 
With the above notation we obtain \\
\\
a) $||u||_{H_{x}^{s}(\mathbb{R})}  =
 ||<\xi>^{s}\widehat{u}\,||_{L^{2}(\mathbb{R})}.$ \\
b) $||u||_{L_{t}^{2}(\mathbb{R}:\,H_{x}^{r}(\mathbb{R}))} =
||<\xi>^{r}\widehat{u}\,||_{L^{2}(\mathbb{R}^{2})}.$ \\
c) $||<D_{x}>^{s}u||_{L^{2}(\mathbb{R})} =
||u||_{H^{s}(\mathbb{R})}.$ \\
d) $||<D_{t}>^{b}<D_{x}>^{r}u||_{L_{x,\,t}^{2}(\mathbb{R}^{2})} =
||u||_{H_{t}^{b}(\mathbb{R}:\,H_{x}^{r}(\mathbb{R}))}.$ \\
e)
$||<D_{x,\,t}>^{s}u||_{L_{t}^{2}(\mathbb{R}:\,H_{x}^{r}(\mathbb{R}))}
= ||<\xi>^{r}\,<|\xi| + |\tau|>^{s}\widehat{u}(\xi,\,\tau)
||_{L^{2}(\mathbb{R}^{2})}.$ 
\end{remark}
We consider the following operators: $L=\partial_{t} +
\partial_{x}^{3}$ and $J=x -
3\,t\,\partial_{x}^{2}$ then $[L,\,J]\equiv L\,J - J\,L=0.$ We
introduce the {\it "generator of dilation"} $P=3\,t\,\partial_{t} +
x\,\partial_{x}$ for the linear part of the coupled system
\eqref{e207}-\eqref{e209} and the {\it "localized dilation
operator"} $P_{0}=3\,t_{0}\,\partial_{t} + x_{0}\,\partial_{x}.$\\
By employing a localization argument, we look at the operator $P$ as a
vector field $P_{0}=3\,t_{0}\,\partial_{t} + x_{0}\,\partial_{x}$
near a fixed point $(x_{0},\,t_{0})\in\mathbb{R}\times\{(-T,\,0)\cup
(0,\,T)\}.$ Since $P_{0}$ is a directional derivative toward to
$(x_{0},\,t),$ we introduce another operator ${\cal
L}_{0}^{3}=t_{0}\,\partial_{x}^{3}$ which plays the role of a
non-tangential vector field to $P_{0}.$ Since $P_{0}$ and ${\cal
L}_{0}$ are linearly independent, the space and time derivative can
be covered by those operator. The main reason why we choose ${\cal
L}_{0}$ is because the corresponding variable coefficients operator
${\cal L}^{3}=t\,\partial_{x}^{3}$ can be treated via the equations
\eqref{e111}-\eqref{e112} and a cut-off procedure enables us to
handle the right hand side of
those.
\begin{remark} For $L$ and $P$ we have the following properties:\\
\\
a) $[L,\,P]\equiv L\,P = (P + 3)L.$ \\
b) $L\,P^{k} = (P + 3)^{k}L.$ \\
c) $(P + 3)^{k}\partial_{x} = \partial_{x}(P + 2)^{k}.$ \\
d) $(P + 3)^{k}\partial_{x}^{3} =
\partial_{x}^{3}P^{k}.$ \\
e) $P_{0}\,P = P\,P_{0} + 3\,P_{0} - 2\,x_{0}\,\partial_{x}.$
\end{remark}
\begin{notation}
The summation $\displaystyle \quad \sum_{\stackrel{k=k_{1} + k_{2} +
k_{3}}{0\leq k_{1},\,k_{2},\,k_{3}\leq k}}$ is simply
abbreviated by $\displaystyle \sum_{k=k_{1} + k_{2} + k_{3}}$.
\end{notation}
Let $P^{k}u=u_{k},$ then
\begin{eqnarray*}
\partial_{t}(P^{k}u) + \partial_{x}^{3}(P^{k}u)
& = & L\,P^{k}u
 =  (P + 3)^{k}Lu
 =  (P + 3)^{k}(\partial_{t}u + \partial_{x}^{3}u) \nonumber \\
& = & -(P + 3)^{k}\left[\frac{a}{2}\;\partial_{x}(u^{2}) +
\frac{b}{2}\;\partial_{x}(v^{2}) + c\,\partial_{x}(u\,v)\right]
\nonumber \\
& = & -\,\frac{a}{2}\,(P + 3)^{k}\partial_{x}(u^{2}) -
\frac{b}{2}\,(P +
3)^{k}\partial_{x}(v^{2}) - c\,(P + 3)^{k}\partial_{x}(u\,v)\nonumber \\
& = & -\,\frac{a}{2}\;\partial_{x}(P + 2)^{k}(u^{2}) -
\frac{b}{2}\;\partial_{x}(P + 2)^{k}(v^{2}) - c\;\partial_{x}(P +
2)^{k}(u\,v).
\end{eqnarray*}
Noting that
$\displaystyle 
(P + 2)^{k}u = \sum_{j=0}^{k}{k\choose j}2^{k - j}P^{j}u.
$
Hence
\begin{eqnarray}
B_{k}^{1}(u,\,u) & = & -\,\frac{a}{2}\;\partial_{x}(P + 2)^{k}(u^{2})  \nonumber \\
& = & -\,\frac{a}{2}\;\partial_{x}\sum_{m=0}^{k}{k\choose m}\;(P +
2)^{m}u\cdot
P^{k - m}u \nonumber \\
& = & -\,\frac{a}{2}\;\partial_{x}\sum_{m=0}^{k}
\sum_{j=0}^{m}{k\choose m}{m\choose j}\;2^{m - j}\;P^{j}u\cdot
P^{k - m}u \nonumber \\
& = & -\,\frac{a}{2}\;\partial_{x}\sum_{m=0}^{k}
\sum_{j=0}^{m}\frac{k!}{(m - j)!\,j!\,(k - m)!}\;2^{m -
j}\;P^{j}u\cdot
P^{k - m}u \nonumber \\
\label{e211}& = & -\,\frac{a}{2}\sum_{k=k_{1} + k_{2} + k_{3}}
\frac{k!}{k_{1}!\;k_{2}!\;k_{3}!}\;2^{k_{1}}\;\partial_{x}\left(u_{k_{2}}\cdot
u_{k_{3}}\right).
\end{eqnarray}
In a similar way
\begin{eqnarray}
\label{e212}B_{k}^{2}(v,\,v) = -\,\frac{b}{2}\;\partial_{x}(P +
2)^{k}(v^{2})
 =  -\,\frac{b}{2}\sum_{k=k_{1}' + k_{2}' +
k_{3}'}
\frac{k!}{k_{1}'!\;k_{2}'!\;k_{3}'!}\;2^{k_{1}'}\;\partial_{x}\left(v_{k_{2}'}\cdot
v_{k_{3}'}\right).
\end{eqnarray}
\begin{eqnarray}
\label{e213}B_{k}^{3}(u,\,v) = c\;\partial_{x}(P + 2)^{k}(u\,v)
 =  -\;c\sum_{k=k_{1}'' + k_{2}'' +
k_{3}''}
\frac{k!}{k_{1}''!\;k_{2}''!\;k_{3}''!}\;2^{k_{1}''}\;\partial_{x}\left(u_{k_{2}''}\cdot
v_{k_{3}''}\right).
\end{eqnarray}
Therefore
\begin{eqnarray}
\lefteqn{\partial_{t}(P^{k}u) + \partial_{x}^{3}(P^{k}u)} \nonumber \\
& = &  -\;\frac{a}{2}\sum_{k=k_{1} + k_{2} + k_{3}}
\frac{k!}{k_{1}!\;k_{2}!\;k_{3}!}\;2^{k_{1}}\;\partial_{x}\left(u_{k_{2}}\cdot
u_{k_{3}}\right) - \frac{b}{2}\sum_{k=k_{1}' + k_{2}' + k_{3}'}
\frac{k!}{k_{1}'!\;k_{2}'!\;k_{3}'!}\;2^{k_{1}'}\;\partial_{x}\left(v_{k_{2}'}\cdot
v_{k_{3}'}\right)\nonumber \\
&  & -\;c\sum_{k=k_{1}'' + k_{2}'' + k_{3}''}
\frac{k!}{k_{1}''!\;k_{2}''!\;k_{3}''!}\;2^{k_{1}''}\;\partial_{x}\left(u_{k_{2}''}\cdot
v_{k_{3}''}\right) \nonumber \\
\label{e214}& = & B_{k}^{1}(u,\,u) + B_{k}^{2}(v,\,v) +
B_{k}^{3}(u,\,v).
\end{eqnarray}
Performing similar calculations as above we obtain
\begin{eqnarray}
\lefteqn{\partial_{t}(P^{k}v) + \partial_{x}^{3}(P^{k}v) }
\nonumber
\\
& = &  -\;\frac{\widetilde{a}}{2}\sum_{k=k_{1} + k_{2} + k_{3}}
\frac{k!}{k_{1}!\;k_{2}!\;
k_{3}!}\;2^{k_{1}}\;\partial_{x}\left(u_{k_{2}}\cdot
u_{k_{3}}\right) - \frac{\widetilde{b}}{2}\sum_{k=k_{1}' + k_{2}' +
k_{3}'} \frac{k!}{k_{1}'!\;k_{2}'!\;
k_{3}'!}\;2^{k_{1}'}\;\partial_{x}\left(v_{k_{2}'}\cdot
v_{k_{3}'}\right) \nonumber \\
&  & -\;\widetilde{c}\sum_{k=k_{1}'' + k_{2}'' + k_{3}''}
\frac{k!}{k_{1}''!\;k_{2}''!\;k_{3}''!}
\;2^{k_{1}''}\;\partial_{x}\left(u_{k_{2}''}\cdot
v_{k_{3}''}\right)\nonumber \\
\label{e215}& = & C_{k}^{1}(u,\,u) + C_{k}^{2}(v,\,v) +
C_{k}^{3}(u,\,v).
\end{eqnarray}
The above nonlinear terms maintain the bilinear structure like that
of the original coupled system of equations of KdV type, since
Leibniz's rule can be applied for operations of $P.$ Now, each
$u_{k}$ and $v_{k}$ satisfies the following system of equations
\begin{eqnarray}
\label{e216}\partial_{t}u_{k} + \partial_{x}^{3}u_{k} & = &
B_{k}^{1}(u,\,u) + B_{k}^{2}(v,\,v) + B_{k}^{3}(u,\,v) \equiv B_{k}\\
\label{e217}\partial_{t}v_{k} + \partial_{x}^{3}v_{k} & = &
C_{k}^{1}(u,\,u) + C_{k}^{2}(v,\,v) +
C_{k}^{3}(u,\,v)\equiv C_{k}\\
\label{e218}u_{k}(x,\,0) = (x\,\partial_{x})^{k}u_{0}(x)\equiv
u_{0}^{k}(x), &  & v_{k}(x,\,0) =
(x\,\partial_{x})^{k}v_{0}(x)\equiv v_{0}^{k}(x).
\end{eqnarray}
In order to obtain a well-posedness result for the system
\eqref{e216}-\eqref{e218} we use Duhamel's principle and we study
the following system of integral equations equivalent to the system
\eqref{e216}-\eqref{e218}
\begin{eqnarray}
\label{e219}\psi(t)\,u_{k} & = & \psi(t)\,V(t)\,u_{0}^{k} -
\psi(t)\int_{0}^{t}V(t - t')\,\psi_{T}(t')\,B_{k}(t')\,dt' \\
\label{e220}\psi(t)\,v_{k} & = & \psi(t)\,V(t)\,v_{0}^{k} -
\psi(t)\int_{0}^{t}V(t - t')\,\psi_{T}(t')\,C_{k}(t')\,dt'
\end{eqnarray}
where $V(t)=e^{-t\,\partial_{x}^{3}}$ is the unitary group
associated with the linear problem and $\psi(t)\in
C_{0}^{\infty}(\mathbb{R}),$ $0\leq\psi\leq 1$ is a cut-off
function such that
\[\psi(t) =\left \{ \begin{array}{ll}
1,& \mbox{if} \quad |t|<1\\
0, & \mbox{if} \quad |t|>2
\end{array}
\qquad \mbox{and} \qquad  \psi_{T}(t)=\psi(t/T)\right. \]
The following results are going to be used several times in the rest of this paper.
\begin{lemma}[\cite{ke0}]. {\it Let $s\in\mathbb{R},$ $a,\,a'\in
(0,\,1/2),$ $b\in (1/2,\,1)$ and $\delta<1.$ Then for any
$k=0,\,1,\,2,\ldots,$ we have}
\begin{eqnarray}
\label{e221}||\psi_{\delta}\phi_{k}||_{X_{-a}^{s}}\leq c\,\delta^{(a
- a')/4(1 -
a')}\,||\phi_{k}||_{X_{-a'}^{s}},\\
\label{e222}||\psi_{\delta}\,V(t)\,\phi_{k}||_{X_{b}^{s}}\leq
c\,\delta^{1/2 - b}||\phi_{k}||_{H^{s}(\mathbb{R})},\\
\label{e223}\left|\left|\psi_{\delta}\int_{0}^{t}V(t -
t')\,F_{k}(t')\,dt'\right|\right|_{X_{b}^{s}}\leq c\,\delta^{1/2 -
b}\,||F_{k}||_{X_{b - 1}^{s}}.
\end{eqnarray}
\end{lemma}
\begin{lemma}[\cite{ke0}]. {\it Let $s>-3/4,$ $b,\,b'\in
(1/2,\,7/12)$ with $b<b'.$ Then for any $k,\,l=0,\,1,\,2,\ldots $ we
have}
\begin{eqnarray}
\label{e224}||\partial_{x}(u_{k}\,v_{l})||_{X_{b' - 1}^{s}}\leq
c\,||v_{k}||_{X_{b}^{s}}\,||v_{l}||_{X_{b}^{s}}.
\end{eqnarray}
\end{lemma}
\begin{lemma}[\cite{ka3}]. {\it Let $s<0,$ $b\in (1/2,\,7/12)$ and
$\psi=\psi(x,\,t)$ be a smooth cut-off function such that the
support of $\psi$ is in $\mathbb{B}_{2}(0)$ and $\psi=1$ on
$\mathbb{B}_{1}(0).$ We set $\psi_{\epsilon}=\psi((x -
x_{0})/\epsilon,\,(t - t_{0})/\epsilon).$ Then for $f\in X_{b}^{s},$
we have}
\begin{eqnarray}
\label{e225}||\psi_{\epsilon}\,f||_{X_{b}^{s}}\leq c\,\epsilon^{-|s|
- 5|b|}||\psi_{\epsilon}||_{X_{|b|}^{|s| + 2\,|b|}}\,||f||_{X_{b}^{s
+ 2\,|b|}},
\end{eqnarray}
where the constant $c$ is independent of $\epsilon$ and $f.$
\end{lemma}
\begin{lemma}[\cite{ka3}]. {\it Let $P$ be the generator of the
dilation and $D_{x,\,t}$ be an operator defined by ${\cal
F}_{\xi,\,\tau}^{-1}<|\tau| + |\xi|>{\cal F}_{x,\,t}.$ We fix an
arbitrary point $(x_{0},\,t_{0})\in \mathbb{R}\times\{(-T,\,0)\cup
(0,\,T)\}.$ Then}\\
1) {\it Suppose that $b\in (0,\,1],$ $r\in (-\infty,\,0]$ and $g\in
X_{b - 1}^{r}$ with $supp\,g\subset
\mathbb{B}_{2\epsilon}(x_{0},\,t_{0})$ and $t\partial_{x}^{3}g,$
$P^{3}g\in X_{b - 1}^{r}.$ If $\epsilon >0$ is sufficiently small,
then we have}
\begin{eqnarray}
\label{e226}||<D_{x,\,t}>^{3b}g||_{L^{2}(\mathbb{R}:\,H^{r}(\mathbb{R}))}\leq
c\,\left(||g||_{X_{b - 1}^{r}} + ||t\partial_{x}^{3}g||_{X_{b -
1}^{r}} + ||P^{3}g||_{X_{b - 1}^{r}}\right)
\end{eqnarray}
{\it where the constant $c=c(x_{0},\,t_{0},\,\epsilon).$}\\
2) {\it If $g\in H^{\mu - 3}(\mathbb{R}^{2})$ with $supp\,g\subset
\mathbb{B}_{2\epsilon}(x_{0},\,t_{0})$ and $t\partial_{x}^{3}g,$
$P^{3}g\in H^{\mu - 3}(\mathbb{R}^{2}).$ Then for small $\epsilon,$
we have}
\begin{eqnarray}
\label{e227}||<D_{x,\,t}>^{\mu}g||_{L^{2}(\mathbb{R}^{2})}\leq
c\,\left(||g||_{H^{\mu - 3}(\mathbb{R}^{2})} +
||t\partial_{x}^{3}g||_{H^{\mu - 3}(\mathbb{R}^{2})} +
||P^{3}g||_{H^{\mu - 3}(\mathbb{R}^{2})}\right)
\end{eqnarray}
{\it where the constant $c=c(x_{0},\,t_{0},\,\epsilon).$}
\end{lemma}
\begin{lemma}[\cite{ka3}]. {\it Let $0\leq s,\,r\leq n/2$ with
$n/2\leq s + r$ and suppose that $f\in H^{s}(\mathbb{R}^{n})$ and $g\in
H^{r}(\mathbb{R}^{n}).$ Then for any $\sigma< s + r - n/2,$ we have
$f\,g\in H^{\sigma}(\mathbb{R}^{n})$ and}
\begin{eqnarray}
\label{e228}||f\,g||_{H^{\sigma}(\mathbb{R}^{n})}\leq
c(\epsilon)\,||f||_{H^{s}(\mathbb{R}^{n})}||g||_{H^{r}(\mathbb{R}^{n})},
\end{eqnarray}
{\it where $\epsilon=s + r - n/2 - \sigma.$}
\end{lemma}
\begin{corollary}[\cite{ka3}]. {\it For $1/2<b<1$ and $-3/4<s<0,$
we have}
\begin{eqnarray}
\label{e229}||\psi\,f||_{X_{b - 1}^{s - 1}}\leq c\,||f||_{X_{b -
1}^{s}}
\end{eqnarray}
where $\psi\in C_{0}^{\infty}(\mathbb{R}^{2})$ and $c$ is
independent of $f.$\end{corollary}
\begin{lemma}[\cite{ka3}]. {\it Let $\psi(x)$ be a smooth cut-off
function in $C_{0}^{\infty}((-2,\,2))$ with $\psi(x)=1$ on
$(-1,\,1).$ We set $\psi_{\epsilon}=\psi(x/\epsilon)$ for
$0<\epsilon <1.$ Then for $r\leq 0,$ and $f\in H^{r},$ we have}
\[||\psi_{\epsilon}\,f||_{H^{r}(\mathbb{R})} \leq\left \{ \begin{array}{ll}
c\,\epsilon^{-\delta}||f||_{H^{r}(\mathbb{R})} &
\quad \mbox{if}\quad -1/2\leq r\leq 0\\
c\,\epsilon^{1/2 + r}||f||_{H^{r}(\mathbb{R})} & \quad
\mbox{if}\quad r<-1/2
\end{array}
\right. \] {\it where $\delta>0$ is an arbitrary small constant and
$c$
is independent of $\epsilon.$}
\end{lemma}
Throughout this paper $c$ is a generic constant, not necessarily the
same at each occasion (it will change from line to line), which
depends in an increasing way on the indicated quantities.
\renewcommand{\theequation}{\thesection.\arabic{equation}}
\setcounter{equation}{0}\section{Existence and Well-Posedness} We
firstly solve the following (slightly general) system of
equations
\begin{eqnarray}
\label{e301}\partial_{t}u_{k} + \partial_{x}^{3}u_{k} & = &
B_{k}^{1}(u,\,u) + B_{k}^{2}(v,\,v) + B_{k}^{3}(u,\,v) \equiv
B_{k}\\
\label{e302}\partial_{t}v_{k} + \partial_{x}^{3}v_{k} & = &
C_{k}^{1}(u,\,u) + C_{k}^{2}(v,\,v) + C_{k}^{3}(u,\,v)\equiv C_{k}\\
\label{e303}u_{k}(x,\,0) = (x\,\partial_{x})^{k}u_{0}(x)\equiv
u_{0}^{k}(x) & , & v_{k}(x,\,0) =
(x\,\partial_{x})^{k}v_{0}(x)\equiv v_{0}^{k}(x)
\end{eqnarray}
where $B_{k}$ and $C_{k}$ are as above.
\begin{definition} 
Let $f=(f_{0},\,f_{1},\,\ldots\,,\,f_{k})$
denotes the infinity series of distributions and define
\begin{eqnarray*}
{\cal A}_{A_{0}}(X_{b}^{s}) \equiv
\left\{f=(f_{0},\,f_{1},\,\ldots\,,\,f_{k}),\,f_{i}\in
X_{b}^{s},\;(i=0,\,1,\,2\ldots)\;\mbox{such that}\;||f||_{{\cal
A}_{A_{0}}(X_{b}^{s})}<+\infty\right\}
\end{eqnarray*}
where
\begin{eqnarray*}
||f||_{{\cal
A}_{A_{0}}(X_{b}^{s})}\equiv\sum_{k=0}^{\infty}\frac{A_{0}^{k}}{k!}\;||f_{k}||_{X_{b}^{s}}.
\end{eqnarray*}
Similarly, for
$u_{0}=\{u_{0}^{0},\,u_{0}^{1},\,\ldots\,,\,u_{0}^{k},\,\ldots\,\}$
and
$v_{0}=\{v_{0}^{0},\,v_{0}^{1},\,\ldots\,,\,v_{0}^{k},\,\ldots\,\}$
we set
\begin{eqnarray*}
||u_{0}||_{{\cal A}_{A_{0}}(H^{s}(\mathbb{R}))}\equiv
\sum_{k=0}^{\infty}\frac{A_{0}^{k}}{k!}\;||u_{0}^{k}||_{H^{s}(\mathbb{R})}
\quad\mbox{and}\quad ||v_{0}||_{{\cal
A}_{A_{0}}(H^{s}(\mathbb{R}))}\equiv
\sum_{k=0}^{\infty}\frac{A_{0}^{k}}{k!}\;||v_{0}^{k}||_{H^{s}(\mathbb{R})}
\end{eqnarray*}
respectively. 
\end{definition}
\begin{remark}
Each solution of the coupled system of
Korteweg de Vries equations is accompanied by the following estimate
\begin{eqnarray*}
||P^{k}u||_{X_{b}^{s}}\leq c\,A_{0}^{k}\,k!,\qquad
\hbox{
and
}\qquad
||P^{k}v||_{X_{b}^{s}}\leq c\,A_{1}^{k}\,k!,\qquad
k=0,\,1,\,2,\ldots
\end{eqnarray*}
\end{remark}
\begin{theorem} {\it Let $-3/4<s,$ $b\in (1/2,\,7/12).$ Suppose
that $u_{0}^{k},\,v_{0}^{k}\in
H^{s}(\mathbb{R})$($k=0,\,1,\,2,\,\ldots$) and satisfies
\begin{eqnarray*}
||u_{0}||_{{\cal A}_{A_{0}}(X_{b}^{s})} =
\sum_{k=0}^{\infty}\frac{A_{0}^{k}}{k!}\;||u_{0}^{k}||_{H^{s}(\mathbb{R})}<+\infty
\quad\mbox{and}\quad ||v_{0}||_{{\cal A}_{A_{0}}(X_{b}^{s})} =
\sum_{k=0}^{\infty}\frac{A_{0}^{k}}{k!}\;||v_{0}^{k}||_{H^{s}(\mathbb{R})}<+\infty.
\end{eqnarray*}
Then there exist
$T=T(||u_{0}^{k}||_{H^{s}(\mathbb{R})},\,||v_{0}^{k}||_{H^{s}(\mathbb{R})})$
and a unique solution $u=(u_{0},\,u_{1},\,\ldots)$ and
$v=v(v_{0},\,v_{1},\,\ldots)$ of the system
\eqref{e301}-\eqref{e303} with $u_{k},\,v_{k}\in
C((-T,\,T):\;H^{s}(\mathbb{R}))\cap X_{b}^{s}$ and
\begin{eqnarray*}
\sum_{k=0}^{\infty}\frac{A_{0}^{k}}{k!}\;||u_{k}||_{X_{b}^{s}(\mathbb{R})}<+\infty,\qquad
\sum_{k=0}^{\infty}\frac{A_{0}^{k}}{k!}\;||v_{k}||_{X_{b}^{s}(\mathbb{R})}<+\infty.
\end{eqnarray*}
Moreover, the map $ (u_{0}^{k},\,v_{0}^{k})\rightarrow
(u(t),\,v(t))$ is Lipschitz continuous, i. e.,
\begin{eqnarray*}
||u(t) - \widetilde{u}(t)||_{{\cal A}_{A_{0}}(X_{b}^{s})} + ||u(t) -
\widetilde{u}(t)||_{C((-T,\,T):\;H^{s}(\mathbb{R}))}\leq
c(T)\,||u_{0} - \widetilde{u}_{0}||_{{\cal
A}_{A_{0}}(H^{s}(\mathbb{R}))}
\end{eqnarray*}
and}
\begin{eqnarray*}
||v(t) - \widetilde{v}(t)||_{{\cal A}_{A_{0}}(X_{b}^{s})} + ||v(t) -
\widetilde{v}(t)||_{C((-T,\,T):\;H^{s}(\mathbb{R}))}\leq
c(T)\,||v_{0} - \widetilde{v}_{0}||_{{\cal
A}_{A_{0}}(H^{s}(\mathbb{R}))}.
\end{eqnarray*}
\end{theorem}
\noindent
{\it Proof.} For given $(u_{0},\,v_{0})\in {\cal
A}_{A_{0}}(H^{s}(\mathbb{R}))\times {\cal
A}_{A_{0}}(H^{s}(\mathbb{R}))$ and $b>1/2,$ let us define,
\begin{eqnarray*}
\mathbb{H}_{R_{1},\,R_{2}} =  \left\{(u,\,v)\in {\cal
A}_{A_{0}}(X_{b}^{s})\times {\cal
A}_{A_{0}}(X_{b}^{s}):\quad||u||_{{\cal A}_{A_{0}}(X_{b}^{s})}\leq
R_{1},\quad ||v||_{{\cal A}_{A_{0}}(X_{b}^{s})}\leq R_{2}\right\}
\end{eqnarray*}
where $R_{1}=2\,c_{0}\,||u_{0}||_{{\cal
A}_{A_{0}}(H^{s}(\mathbb{R}))}$ and
$R_{2}=2\,c_{0}\,||v_{0}||_{{\cal A}_{A_{0}}(H^{s}(\mathbb{R}))}.$
Then $\mathbb{H}_{R_{1},\,R_{2}}$ is a complete metric space with
norm
\begin{eqnarray*}
||(u,\,v)||_{\mathbb{H}_{R_{1},\,R_{2}}}=||u||_{{\cal
A}_{A_{0}}(X_{b}^{s})} + ||v||_{{\cal A}_{A_{0}}(X_{b}^{s})}.
\end{eqnarray*}
Without loss of generality, we may assume that that $R_{1}>1$ and $R_{2}>1.$
For $(u,\,v)\in \mathbb{H}_{R_{1},\,R_{2}},$ let us define the
maps,
\begin{eqnarray}
\label{e126}\Phi_{u_{0}}^{k}(u,\,v) & = & \psi(t)\,V(t)\,u_{0}^{k}
-
\psi(t)\int_{0}^{t}V(t - t')\,\psi_{T}(t')\,B_{k}(t')\,dt' \\
\label{e127}\Psi_{v_{0}}^{k}(u,\,v) & = & \psi(t)\,V(t)\,v_{0}^{k}
- \psi(t)\int_{0}^{t}V(t - t')\,\psi_{T}(t')\,C_{k}(t')\,dt'.
\end{eqnarray}
We prove that $\Phi\times\Psi$ maps $\mathbb{H}_{R_{1},\,R_{2}}$
into $\mathbb{H}_{R_{1},\,R_{2}}$ and it is a contraction. In fact,
using lemma 2.1 and lemma 2.2 we have
\begin{eqnarray*}
||\Phi_{u_{0}}^{k}(u,\,v)||_{X_{b}^{s}} & = &
||\psi(t)\,V(t)\,u_{0}^{k}||_{X_{b}^{s}} +
\left|\left|\psi(t)\int_{0}^{t}V(t -
t')\,\psi_{T}(t')\,B_{k}(t')\,dt'\right|\right|_{X_{b}^{s}}\\
& \leq & c_{0}\,||u_{0}^{k}||_{H^{s}(\mathbb{R})} +
c_{1}\,T^{\mu}\,||B_{k}||_{X_{b' - 1}^{s}}\\
& \leq & c_{0}\,||u_{0}^{k}||_{H^{s}(\mathbb{R})} +
c_{1}\,T^{\mu}\,\frac{a}{2}\sum_{k=k_{1} + k_{2} +
k_{3}}\frac{k!}{k_{1}!\,k_{2}!\,k_{3}!}\;2^{k_{1}}\,
||u_{k_{2}}||_{X_{b}^{s}}\,||u_{k_{3}}||_{X_{b}^{s}}\\
&  & +\;c_{1}\,T^{\mu}\,\frac{b}{2}\sum_{k=k_{1}' + k_{2}' +
k_{3}'}\frac{k!}{k_{1}'!\,k_{2}'!\,k_{3}'!}\;2^{k_{1}'}\,
||v_{k_{2}'}||_{X_{b}^{s}}\,||v_{k_{3}'}||_{X_{b}^{s}}\\
&  & +\;c_{1}\,T^{\mu}\,c\sum_{k=k_{1}'' + k_{2}'' +
k_{3}''}\frac{k!}{k_{1}''!\,k_{2}''!\,k_{3}''!}\;2^{k_{1}''}\,
||u_{k_{2}''}||_{X_{b}^{s}}\,||v_{k_{3}''}||_{X_{b}^{s}}.
\end{eqnarray*}
Applying a sum over $k$ we have
\begin{eqnarray*}
\lefteqn{\sum_{k=0}^{\infty}\frac{A_{0}^{k}}{k!}\;||\Phi_{u_{0}}^{k}(u,\,v)||_{X_{b}^{s}}}\\
& \leq &
c_{0}\,\sum_{k=0}^{\infty}\frac{A_{0}^{k}}{k!}\;||u_{0}^{k}||_{H^{s}(\mathbb{R})}
+
c_{1}\,T^{\mu}\,\frac{a}{2}\sum_{k=0}^{\infty}\frac{A_{0}^{k}}{k!}\sum_{k=k_{1}
+ k_{2} + k_{3}}\frac{k!}{k_{1}!\,k_{2}!\,k_{3}!}\;2^{k_{1}}\,
||u_{k_{2}}||_{X_{b}^{s}}\,||u_{k_{3}}||_{X_{b}^{s}}\\
&  &
+\;c_{1}\,T^{\mu}\,\frac{b}{2}\sum_{k=0}^{\infty}\frac{A_{0}^{k}}{k!}\sum_{k=k_{1}'
+ k_{2}' + k_{3}'}\frac{k!}{k_{1}'!\,k_{2}'!\,k_{3}'!}\;2^{k_{1}'}\,
||v_{k_{2}'}||_{X_{b}^{s}}\,||v_{k_{3}'}||_{X_{b}^{s}}\\
&  &
+\;c_{1}\,T^{\mu}\,c\sum_{k=0}^{\infty}\frac{A_{0}^{k}}{k!}\sum_{k=k_{1}''
+ k_{2}'' +
k_{3}''}\frac{k!}{k_{1}''!\,k_{2}''!\,k_{3}''!}\;2^{k_{1}''}\,
||u_{k_{2}''}||_{X_{b}^{s}}\,||v_{k_{3}''}||_{X_{b}^{s}}\\
& \leq & c_{0}\;||u_{0}||_{{\cal A}_{A_{0}}(H^{s}(\mathbb{R}))} +
c_{1}\,T^{\mu}\,\frac{a}{2}\sum_{k=0}^{\infty}\sum_{k=k_{1} + k_{2}
+
k_{3}}\,2^{k_{1}}\;\frac{A_{0}^{k_{1}}}{k_{1}!}\;\frac{A_{0}^{k_{2}}}{k_{2}!}\;
||u_{k_{2}}||_{X_{b}^{s}}\;\frac{A_{0}^{k_{3}}}{k_{3}!}\;||u_{k_{3}}||_{X_{b}^{s}}\\
&  & +\;c_{1}\,T^{\mu}\,\frac{b}{2}\sum_{k=0}^{\infty}\sum_{k=k_{1}'
+ k_{2}' +
k_{3}'}\,2^{k_{1}'}\;\frac{A_{0}^{k_{1}'}}{k_{1}'!}\;\frac{A_{0}^{k_{2}'}}{k_{2}'!}\;
||v_{k_{2}'}||_{X_{b}^{s}}\;\frac{A_{0}^{k_{3}'}}{k_{3}'!}\;||v_{k_{3}'}||_{X_{b}^{s}}\\
&  & +\;c_{1}\,T^{\mu}\,c\sum_{k=0}^{\infty}\sum_{k=k_{1}'' +
k_{2}'' +
k_{3}''}\,2^{k_{1}''}\;\frac{A_{0}^{k_{1}''}}{k_{1}'!}\;\frac{A_{0}^{k_{2}''}}{k_{2}''!}\;
||u_{k_{2}''}||_{X_{b}^{s}}\;\frac{A_{0}^{k_{3}''}}{k_{3}''!}\;||v_{k_{3}''}||_{X_{b}^{s}}\\
& \leq & c_{0}\;||u_{0}||_{{\cal A}_{A_{0}}(H^{s}(\mathbb{R}))} +
c_{1}\,T^{\mu}\,\frac{a}{2}\sum_{k_{1}=0}^{\infty}2^{k_{1}}\;\frac{A_{0}^{k_{1}}}{k_{1}!}
\sum_{k_{2}=0}^{\infty}\;\frac{A_{0}^{k_{2}}}{k_{2}!}\;||u_{k_{2}}||_{X_{b}^{s}}
\sum_{k_{3}=0}^{\infty}\;\frac{A_{0}^{k_{3}}}{k_{3}!}\;||u_{k_{3}}||_{X_{b}^{s}}\\
&  &
+\;c_{1}\,T^{\mu}\,\frac{b}{2}\sum_{k_{1}'=0}^{\infty}2^{k_{1}'}\;\frac{A_{0}^{k_{1}'}}{k_{1}'!}
\sum_{k_{2}'=0}^{\infty}\;\frac{A_{0}^{k_{2}'}}{k_{2}'!}\;||v_{k_{2}'}||_{X_{b}^{s}}
\sum_{k_{3}'=0}^{\infty}\;\frac{A_{0}^{k_{3}'}}{k_{3}'!}\;||v_{k_{3}'}||_{X_{b}^{s}}\\
&  &
+\;c_{1}\,T^{\mu}\,c\sum_{k_{1}''=0}^{\infty}2^{k_{1}''}\;\frac{A_{0}^{k_{1}''}}{k_{1}''!}
\sum_{k_{2}''=0}^{\infty}\;\frac{A_{0}^{k_{2}''}}{k_{2}''!}\;||u_{k_{2}''}||_{X_{b}^{s}}
\sum_{k_{3}''=0}^{\infty}\;\frac{A_{0}^{k_{3}''}}{k_{3}''!}\;||v_{k_{3}''}||_{X_{b}^{s}}\\
& = & c_{0}\;||u_{0}||_{{\cal A}_{A_{0}}(H^{s}(\mathbb{R}))} +
c_{1}\,T^{\mu}\,\frac{a}{2}\;e^{2\,A_{0}}\;||u||_{{\cal
A}_{A_{0}}(X_{b}^{s})}^{2}\\
&  & +\;c_{1}\,T^{\mu}\,\frac{b}{2}\;e^{2\,A_{0}}\;||v||_{{\cal
A}_{A_{0}}(X_{b}^{s})}^{2} +
c_{1}\,T^{\mu}\,c\;e^{2\,A_{0}}\;||u||_{{\cal
A}_{A_{0}}(X_{b}^{s})}\;||v||_{{\cal A}_{A_{0}}(X_{b}^{s})}.
\end{eqnarray*}
Hence, choosing $d=\max\{a/2,\,b/2,\,c\}$ we have
\begin{eqnarray}
\lefteqn{||\Phi_{u_{0}}(u,\,v)||_{{\cal
A}_{A_{0}}(X_{b}^{s})} \leq  c_{0}\;||u_{0}||_{{\cal A}_{A_{0}}(H^{s}(\mathbb{R}))}} \nonumber \\
&  & +\;c_{1}\,T^{\mu}\,d\;e^{2\,A_{0}}\;\left[||u||_{{\cal
A}_{A_{0}}(X_{b}^{s})}^{2} + ||v||_{{\cal A}_{A_{0}}(X_{b}^{s})}^{2}
+ ||u||_{{\cal A}_{A_{0}}(X_{b}^{s})}\;||v||_{{\cal
A}_{A_{0}}(X_{b}^{s})}\right] \nonumber \\
\label{e306}& \leq & c_{0}\;||u_{0}||_{{\cal
A}_{A_{0}}(H^{s}(\mathbb{R}))} +
\frac{3}{2}\;c_{1}\,d\;T^{\mu}\;e^{2\,A_{0}}\;\left[||u||_{{\cal
A}_{A_{0}}(X_{b}^{s})}^{2} + ||v||_{{\cal A}_{A_{0}}(X_{b}^{s})}^{2}
\right].
\end{eqnarray}
In a similar way, choosing
$\widetilde{d}=\max\{\widetilde{a}/2,\,\widetilde{b}/2,\,\widetilde{c}\}$
we have
\begin{eqnarray}
\label{e307}||\Psi_{v_{0}}(u,\,v)||_{{\cal A}_{A_{0}}(X_{b}^{s})}
 \leq  c_{0}\;||v_{0}||_{{\cal A}_{A_{0}}(H^{s}(\mathbb{R}))} +
\frac{3}{2}\;c_{2}\,\widetilde{d}\;T^{\mu}\;e^{2\,A_{0}}\;\left[||u||_{{\cal
A}_{A_{0}}(X_{b}^{s})}^{2} + ||v||_{{\cal A}_{A_{0}}(X_{b}^{s})}^{2}
\right].
\end{eqnarray}
If we choose $T$ such that
\begin{eqnarray*}
T^{\mu}\leq \frac{1}{3\;\max\{c_{1},\,c_{2}\}\;(R_{1} + R_{2})^{2}}
\end{eqnarray*}
Then we obtain in \eqref{e306} and \eqref{e307}
\begin{eqnarray*}
||\Phi_{u_{0}}(u,\,v)||_{{\cal A}_{A_{0}}(X_{b}^{s})}\leq
R_{1}\qquad\mbox{and}\qquad ||\Psi_{v_{0}}(u,\,v)||_{{\cal
A}_{A_{0}}(X_{b}^{s})}\leq R_{2}.
\end{eqnarray*}
Therefore, $(\Phi_{u_{0}},\,\Psi_{v_{0}})\in
\mathbb{H}_{R_{1},\,R_{2}}.$ We show that $\Phi_{u_{0}}\times
\Psi_{v_{0}}:(u,\,v)\rightarrow
(\Phi_{u_{0}}(u,\,v),\,\Psi_{v_{0}}(u,\,v))$
is a contraction. \\
Let
$(u,\,v),\;(\widetilde{u},\,\widetilde{v})\in\mathbb{H}_{R_{1},\,R_{2}},$
then as above we get for $d=\max\{a/2,\,b/2,\,c\}$
\begin{eqnarray}
\lefteqn{||\Phi_{u_{0}}(u,\,v) -
\Phi_{u_{0}}(\widetilde{u},\,\widetilde{v})||_{{\cal
A}_{A_{0}}(X_{b}^{s})}} \nonumber \\
\label{e308}& \leq &
\frac{3}{2}\;c_{1}\,d\;T^{\mu}\;e^{2\,A_{0}}\;(R_{1} +
R_{2})\,\left[||u - \widetilde{u}||_{{\cal A}_{A_{0}}(X_{b}^{s})} +
||v - \widetilde{v}||_{{\cal A}_{A_{0}}(X_{b}^{s})}\right].
\end{eqnarray}
In a similar way, choosing
$\widetilde{d}=\max\{\widetilde{a}/2,\,\widetilde{b}/2,\,\widetilde{c}\}$
we have
\begin{eqnarray}
\lefteqn{||\Psi_{v_{0}}(u,\,v) -
\Psi_{v_{0}}(\widetilde{u},\,\widetilde{v})||_{{\cal
A}_{A_{0}}(X_{b}^{s})}} \nonumber \\
\label{e309}& \leq &
\frac{3}{2}\;c_{2}\,\widetilde{d}\;T^{\mu}\;e^{2\,A_{0}}\;(R_{1} +
R_{2})\,\left[||u - \widetilde{u}||_{{\cal A}_{A_{0}}(X_{b}^{s})} +
||v - \widetilde{v}||_{{\cal A}_{A_{0}}(X_{b}^{s})}\right].
\end{eqnarray}
Choosing $T^{\mu}$ small enough, such that
\begin{eqnarray*}
T^{\mu}\leq \frac{1}{6\,\max\{c_{1},\,c_{2}\}\,(R_{1} + R_{2})^{2}}
\end{eqnarray*}
we obtain
\begin{eqnarray}
\label{e310}||\Phi_{u_{0}}(u,\,v) -
\Phi_{u_{0}}(\widetilde{u},\,\widetilde{v})||_{{\cal
A}_{A_{0}}(X_{b}^{s})} \leq \frac{1}{4}\left[||u -
\widetilde{u}||_{{\cal A}_{A_{0}}(X_{b}^{s})} + ||v -
\widetilde{v}||_{{\cal A}_{A_{0}}(X_{b}^{s})}\right].
\end{eqnarray}
In a similar way
\begin{eqnarray}
\label{e311}||\Psi_{v_{0}}(u,\,v) -
\Psi_{v_{0}}(\widetilde{u},\,\widetilde{v})||_{{\cal
A}_{A_{0}}(X_{b}^{s})} \leq  \frac{1}{4}\left[||u -
\widetilde{u}||_{{\cal A}_{A_{0}}(X_{b}^{s})} + ||v -
\widetilde{v}||_{{\cal A}_{A_{0}}(X_{b}^{s})}\right].
\end{eqnarray}
Therefore the map $\Phi_{u_{0}}\times \Psi_{v_{0}}$ is a contraction
and we obtain a unique fixed point $(u,\,v)$ which solves the
initial value problem \eqref{e301}-\eqref{e303} for $T<T^{\mu}.$ The
rest of the proof follows a standard argument.
\begin{corollary} 
{\it Let $-3/4<s,$ $b\in (1/2,\,7/12).$ Suppose
that $(x\,\partial_{x})^{k}u_{0},\,(x\,\partial_{x})^{k}v_{0}\in
H^{s}(\mathbb{R})$($k=0,\,1,\,2,\,\ldots$) and that
\begin{eqnarray*}
\sum_{k=0}^{\infty}\frac{A_{0}^{k}}{k!}\;||u_{0}^{k}||_{H^{s}(\mathbb{R})}<+\infty
\quad\mbox{and}\quad
\sum_{k=0}^{\infty}\frac{A_{0}^{k}}{k!}\;||v_{0}^{k}||_{H^{s}(\mathbb{R})}<+\infty.
\end{eqnarray*}
Then there exist
$T=T(||u_{0}^{k}||_{H^{s}(\mathbb{R})},\,||v_{0}^{k}||_{H^{s}(\mathbb{R})})$
and a unique solution $(u,\,v)$ of the coupled system equations KdV
type \eqref{e101}-\eqref{e103} with $u,\,v\in
C((-T,\,T):\;H^{s}(\mathbb{R}))\cap X_{b}^{s}$ and
\begin{eqnarray*}
\sum_{k=0}^{\infty}\frac{A_{1}^{k}}{k!}\;||P^{k}u||_{X_{b}^{s}(\mathbb{R})}<+\infty,\qquad
\sum_{k=0}^{\infty}\frac{A_{1}^{k}}{k!}\;||P^{k}v||_{X_{b}^{s}(\mathbb{R})}<+\infty.
\end{eqnarray*}
Moreover, the map $ (u_{0},\,v_{0})\rightarrow (u(t),\,v(t))$ is
Lipschitz continuous in the following sense:
\begin{eqnarray*}
||P^{k}u(t) - P^{k}\widetilde{u}(t)||_{X_{b}^{s}} + ||P^{k}u(t) -
P^{k}\widetilde{u}(t)||_{C((-T,\,T):\;H^{s}(\mathbb{R}))}\leq
c(T)\sum_{k=0}^{\infty}\frac{A_{0}^{k}}{k!}\;||(x\,\partial_{x})^{k}(u_{0}
- \widetilde{u}_{0})||_{H^{s}(\mathbb{R})}
\end{eqnarray*}
and}
\begin{eqnarray*}
||v(t) - \widetilde{v}(t)||_{X_{b}^{s}} + ||v(t) -
\widetilde{v}(t)||_{C((-T,\,T):\;H^{s}(\mathbb{R}))}\leq
c(T)\sum_{k=0}^{\infty}\frac{A_{0}^{k}}{k!}\;||(x\,\partial_{x})^{k}(v_{0}
- \widetilde{v}_{0})||_{H^{s}(\mathbb{R})}.
\end{eqnarray*}
\end{corollary}
\renewcommand{\theequation}{\thesection.\arabic{equation}}
\setcounter{equation}{0}\section{The main result} In this section we
prove the analyticity of the solution obtained in the previous section.
We treat the solution $u_{k}\equiv P^{k}u$ and $v_{k}\equiv P^{k}v$
as if they satisfy the coupled system of equations
\eqref{e301}-\eqref{e303} in the classical sense. This can be
justified by a proper approximation procedure. The following results
are going to be used in this section. Let $(x_{0},\,t_{0})$ be
arbitrarily taken in $\mathbb{R}\times \{(-T,\,0)\cup (0,\,T)\}.$ By
$\psi(x,\,t)$ we denote a smooth cut-off function in
$C_{0}^{\infty}(\mathbb{B}_{1}(0))$ and $\psi_{\epsilon}=\psi((x -
x_{0})/\epsilon,\,(t - t_{0})/\epsilon).$ 

Let $\psi$ be a smooth cut-off function around the freezing point
$(x_{0},\,t_{0})$ with $supp\,\psi \subset
C_{0}^{\infty}(\mathbb{B}_{\epsilon}(x_{0},\,t_{0})).$ 
\begin{proposition} {\it For the cut-off function $\psi$ defined
above, there exists a positive constant $c$ and $A$ such that}
\begin{eqnarray}
\label{e401}||\psi\,P^{k}u||_{L_{x,\,t}^{2}(\mathbb{R}^{2})}\leq
c\,A^{k}\,(k!)^{2},\qquad k=0,\,1,\,2,\,\ldots\\
\label{e402}||\psi\,P^{k}v||_{L_{x,\,t}^{2}(\mathbb{R}^{2})}\leq
c\,A^{k}\,(k!)^{2},\qquad k=0,\,1,\,2,\,\ldots
\end{eqnarray}
\end{proposition}
\noindent
{\it Proof.} Using \eqref{e226} with $r=s - 1,$ we obtain
\begin{eqnarray}
\label{e403}||<D_{x,\,t}>^{3b}\psi
P^{k}u||_{L_{t}^{2}(\mathbb{R}:\,H_{x}^{s - 1}(\mathbb{R}))} \leq
c\left(||\psi u_{k}||_{X_{b - 1}^{s - 1}} +
||t\,\partial_{x}^{3}(\psi u_{k})||_{X_{b - 1}^{s - 1}} +
||P^{3}(\psi u_{k})||_{X_{b - 1}^{s - 1}}\right).
\end{eqnarray}
Each term in \eqref{e403} is estimated separately. For the first term in
the right hand side we use Lemma 2.3. Indeed,
\begin{eqnarray}
\label{e404}||\psi\,u_{k}||_{X_{b - 1}^{s - 1}}\leq
||\psi\,u_{k}||_{X_{b - 1}^{s}}\leq c\,||\psi||_{X_{|b - 1|}^{|s| +
2|b - 1|}}||u_{k}||_{X_{b}^{s}}\leq c(\psi)A_{1}^{k}\,k!.\qquad
k=0,\,1,\,2,\,\ldots
\end{eqnarray}
The third term is estimated again using Corollary 2.6.
\begin{eqnarray}
||P^{3}(\psi\,u_{k})||_{X_{b - 1}^{s - 1}} & \leq &
\sum_{l=0}^{3}\frac{3!}{l\,(l - 3)!}\,||(P^{3 - l}
\psi)\,P^{l}u_{k}||_{X_{b - 1}^{s}}\nonumber \\
& \leq & c(\psi)\sum_{l=0}^{3}\frac{3!}{l\,(l - 3)!}\,||P^{l}u_{k}||_{X_{b}^{s}}\nonumber \\
& \leq & c\sum_{l=0}^{3}\frac{3!}{l\,(l - 3)!}\,||P^{k + l}u||_{X_{b}^{s}}\nonumber \\
& = & c\sum_{l=0}^{3}A_{1}^{k + l}\,(k + l)!\nonumber \\
\label{e405}& \leq & c\,A_{2}^{k}\,k!.\qquad k=0,\,1,\,2,\,\ldots
\end{eqnarray}
For the second term, we use \eqref{e301} to reduce the third
derivative in space to the dilation operator $P.$
 Since the generator of dilation is $Pu_{k} =
3\,t\,\partial_{t}u_{k} + x\,\partial_{x}u_{k}$ we obtain
\begin{eqnarray}
\label{e406}t\,\partial_{t}u_{k} = \frac{1}{3}\,Pu_{k} -
\frac{1}{3}\,x\,\partial_{x}u_{k}.
\end{eqnarray}
Multiplying \eqref{e301} by $\psi\,t,$ we have
\begin{eqnarray}
\label{e407}\psi\,t\,\partial_{t}u_{k} +
\psi\,t\,\partial_{x}^{3}u_{k} = \psi\,t\,B_{k}.
\end{eqnarray}
Replacing \eqref{e406} in \eqref{e407} we obtain
\begin{eqnarray}
\label{e408}\psi\,t\,\partial_{x}^{3}u_{k} =
-\,\frac{1}{3}\,\psi\,Pu_{k} +
\frac{1}{3}\,\psi\,x\,\partial_{x}u_{k} + \psi\,t\,B_{k}.
\end{eqnarray}
hence
\begin{eqnarray}
||\psi\,t\,\partial_{x}^{3}u_{k}||_{X_{b - 1}^{s - 1}} & = &
\frac{1}{3}\,||\psi\,Pu_{k}||_{X_{b - 1}^{s - 1}} +
\frac{1}{3}\,||\psi\,x\,\partial_{x}u_{k}||_{X_{b - 1}^{s - 1}} +
||\psi\,t\,B_{k}||_{X_{b - 1}^{s - 1}}\nonumber \\
\label{e409}& = & F_{1} + F_{2} + F_{3}.
\end{eqnarray}
Using the assumption in the Theorem, we have
\begin{eqnarray}
F_{1}  =  \frac{1}{3}\,||\psi\,Pu_{k}||_{X_{b - 1}^{s - 1}} & \leq &
c\,||\psi||_{X_{1 - b}^{-s}}||P^{k + 1}u||_{X_{b - 1}^{s}} \leq
c\,||P^{k + 1}u||_{X_{b}^{s}}\nonumber
\\
& \leq & c\,A_{3}^{k + 1}(k + 1)!
\label{e410} \ \leq  \ c\,A_{4}^{k}\,k!.
\end{eqnarray}
Similarly, we obtain
\begin{eqnarray}
F_{2} = \frac{1}{3}\,||\psi\,x\,\partial_{x}u_{k}||_{X_{b - 1}^{s -
1}} & \leq & \frac{1}{3}\,||\partial_{x}(\psi\,x\,u_{k})||_{X_{b -
1}^{s - 1}} + \frac{1}{3}\,||\partial_{x}(\psi\,x)\,u_{k})||_{X_{b -
1}^{s - 1}}\nonumber \\
& \leq & \frac{1}{3}\,||\partial_{x}(\psi\,x\,v_{k})||_{X_{b -
1}^{s}} + c\,||\partial_{x}(\psi\,x)||_{X_{1 -
b}^{-s}}||u_{k}||_{X_{b -
1}^{s}}\nonumber \\
& \leq & \frac{1}{3}\,||\psi\,x||_{X_{b}^{s}}||u_{k}||_{X_{b}^{s}} +
c\,||\partial_{x}(\psi\,x)||_{X_{1 - b}^{-s}}||u_{k}||_{X_{b -
1}^{s}}\nonumber \\
& \leq & c\,\left(||\psi\,x||_{X_{b}^{s}} +
||\partial_{x}(\psi\,x)||_{X_{1 - b}^{-s}}\right)\,A_{5}^{k}k!
\label{e411}\ \leq \ c\,A_{6}^{k}k!.
\end{eqnarray}
Using Lemma 2.3 and 2.2, we have
\begin{eqnarray*}
F_{3} = ||\psi\,t\,B_{k}||_{X_{b - 1}^{s - 1}} &\leq&
c\,||\psi||_{X_{b - 1}^{-s}}\,||B_{k}^{1} + B_{k}^{2} +
B_{k}^{3}||_{X_{b - 1}^{s}}  \nonumber \\
& \leq & c\,\left(||B_{k}^{1}||_{X_{b - 1}^{s}} +
||B_{k}^{2}||_{X_{b - 1}^{s}} + ||B_{k}^{3}||_{X_{b - 1}^{s}}\right)
\end{eqnarray*}
Then replacing $B_{k}^{1}$, $B_{k}^{2}$ and $B_{k}^{3}$ in \eqref{e211}, \eqref{e212} and \eqref{e213} we deduce
\begin{eqnarray}
F_3 & \leq & c\sum_{k=k_{1} + k_{2} +
k_{3}}\frac{k!}{k_{1}!\,k_{2}!\,k_{3}!}\,2^{k_{1}}\,||u_{k_{2}}||_{X_{b}^{s}}
\,||u_{k_{3}}||_{X_{b}^{s}} + c\sum_{k=k_{1}' + k_{2}' +
k_{3}'}\frac{k!}{k_{1}'!\,k_{2}'!\,k_{3}'!}\,2^{k_{1}'}\,||v_{k_{2}'}||_{X_{b}^{s}}
\,||v_{k_{3}'}||_{X_{b}^{s}} \nonumber \\
&  & +\;c\sum_{k=k_{1}'' + k_{2}'' +
k_{3}''}\frac{k!}{k_{1}''!\,k_{2}''!\,k_{3}''!}\,2^{k_{1}''}\,||u_{k_{2}''}||_{X_{b}^{s}}
\,||v_{k_{3}''}||_{X_{b}^{s}} \nonumber \\
& \leq & c\sum_{k=k_{1} + k_{2} +
k_{3}}\frac{k!}{k_{1}!\,k_{2}!\,k_{3}!}\,2^{k_{1}}\,A_{7}^{k_{2}}\cdot
k_{2}!\;A_{7}^{k_{3}}\cdot k_{3}! + c\sum_{k=k_{1}' + k_{2}' +
k_{3}'}\frac{k!}{k_{1}'!\,k_{2}'!\,k_{3}'!}\,2^{k_{1}'}\,A_{8}^{k_{2}'}\cdot
k_{2}'!\;A_{8}^{k_{3}'}\cdot k_{3}'! \nonumber \\
&  & +\;c\sum_{k=k_{1}'' + k_{2}'' +
k_{3}''}\frac{k!}{k_{1}''!\,k_{2}''!\,k_{3}''!}\,2^{k_{1}''}\,A_{9}^{k_{2}''}\cdot
k_{2}''!\;A_{10}^{k_{3}''}\cdot k_{3}''! \nonumber \\
& \leq & c\sum_{k=k_{1} + k_{2} +
k_{3}}\frac{k!}{k_{1}!}\,2^{k_{1}}\,A_{7}^{k_{2} + k_{3}} +
c\sum_{k=k_{1}' + k_{2}' +
k_{3}'}\frac{k!}{k_{1}'!}\,2^{k_{1}'}\,A_{8}^{k_{2}' + k_{3}'} \nonumber \\
&  & +\;c\sum_{k=k_{1}'' + k_{2}'' +
k_{3}''}\frac{k!}{k_{1}''!}\,2^{k_{1}''}\,A_{9}^{k_{2}''}\cdot
A_{10}^{k_{3}''} \nonumber \\
& \leq & c\,k!\;A_{7}^{k}\sum_{k_{1}=0}^{k}\sum_{k_{2}=0}^{k -
k_{1}}\frac{1}{k_{1}!}\,2^{k_{1}}\,A_{7}^{-\,k_{1}} +
c\,k!\;A_{8}^{k}\sum_{k_{1}'=0}^{k}\sum_{k_{2}'=0}^{k -
k_{1}'}\frac{1}{k_{1}'!}\,2^{k_{1}'}\,A_{8}^{-\,k_{1}'} \nonumber \\
&  & +\;c\,k!\sum_{k=k_{1}'' + k_{2}'' +
k_{3}''}\frac{1}{k_{1}''!}\,2^{k_{1}''}\,A_{9}^{k_{2}''}\cdot
A_{10}^{k_{3}''} \nonumber \\
& \leq & c\,k!\;A_{7}^{k}\sum_{k_{1}=0}^{k}\sum_{k_{2}=0}^{k -
k_{1}}\frac{\left(\frac{2}{A_{7}}\right)^{k_{1}}}{k_{1}!} +
c\,k!\;A_{8}^{k}\sum_{k_{1}'=0}^{k}\sum_{k_{2}'=0}^{k -
k_{1}'}\frac{\left(\frac{2}{A_{8}}\right)^{k_{1}'}}{k_{1}'!} +
c\,k!\sum_{k=k_{1}'' + k_{2}'' +
k_{3}''}\frac{1}{k_{1}''!}\,2^{k_{1}''}\,A_{9}^{k_{2}''}\cdot
A_{10}^{k_{3}''} \nonumber \\
& \leq & c\,e^{2/A_{7}}\;A_{7}^{k}\cdot k! +
c\,e^{2/A_{8}}\;A_{8}^{k}\cdot k! + c\,k!\sum_{k=k_{1}'' + k_{2}'' +
k_{3}''}\frac{1}{k_{1}''!}\,2^{k_{1}''}\,A_{9}^{k_{2}''}\cdot
A_{10}^{k_{3}''} \nonumber \\
\label{e412}& \leq & c\left(e^{2/A_{7}} +
e^{2/A_{8}}\right)A_{11}\cdot k! + c\,k!\sum_{k=k_{1}'' + k_{2}'' +
k_{3}''}\frac{1}{k_{1}''!}\,2^{k_{1}''}\,A_{9}^{k_{2}''}\cdot
A_{10}^{k_{3}''}.\qquad k=0,\,1,\,2,\,\ldots
\end{eqnarray}
Hence, from \eqref{e410}, \eqref{e411} and \eqref{e412} in
\eqref{e409} we obtain that there exists a positive constant $c$ and
$A_{11}$ such that
\begin{eqnarray}
\label{e413}||\psi\,t\,\partial_{x}^{3}u_{k}||_{X_{b - 1}^{s -
1}}\leq c\,A_{11}\cdot k! + c\,k!\sum_{k=k_{1}'' + k_{2}'' +
k_{3}''}\frac{1}{k_{1}''!}\,2^{k_{1}''}\,A_{9}^{k_{2}''}\cdot
A_{10}^{k_{3}''},\qquad k=0,\,1,\,2,\,\ldots
\end{eqnarray}
On the other hand, using  $\partial_{x}^{3}(\psi\cdot f) =
\psi\cdot\partial_{x}^{3}f +
3\,\partial_{x}^{2}(\partial_{x}\psi\cdot f) -
3\,\partial_{x}(\partial_{x}^{2}\psi\cdot f) +
\partial_{x}^{3}\psi\cdot f$
we have that
\begin{eqnarray}
||t\,\partial_{x}^{3}(\psi\cdot u_{k})||_{X_{b - 1}^{s - 1}} & \leq
& ||t\,\psi\cdot\partial_{x}^{3}u_{k}||_{X_{b - 1}^{s - 1}} +
3\,||\partial_{x}^{2}(t\,\partial_{x}\psi\cdot u_{k})||_{X_{b -
1}^{s -
1}} \nonumber \\
\label{e414}&  & +\;3\,||\partial_{x}(t\,\partial_{x}^{2}\psi\cdot
u_{k})||_{X_{b - 1}^{s - 1}} + ||t\,\partial_{x}^{3}\psi\cdot
u_{k}||_{X_{b - 1}^{s - 1}}.
\end{eqnarray}
Using Lemma 2.2 and Lemma 2.3 we obtain
\begin{eqnarray}
||\partial_{x}^{2}(t\,\partial_{x}\psi\cdot u_{k})||_{X_{b - 1}^{s -
1}} & \leq & ||\partial_{x}(t\,\partial_{x}\psi\cdot u_{k})||_{X_{b
- 1}^{s}} \leq
c\,||t\,\partial_{x}\psi||_{X_{b}^{s}}||u_{k}||_{X_{b}^{s}}\nonumber
\\
\label{e415}& \leq & c\,A_{10}^{k}\,k!
\end{eqnarray}
\begin{eqnarray}
||\partial_{x}(t\,\partial_{x}^{2}\psi\cdot u_{k})||_{X_{b - 1}^{s -
1}} & \leq & ||\partial_{x}(t\,\partial_{x}^{2}\psi\cdot
u_{k})||_{X_{b - 1}^{s}} \leq
c\,||t\,\partial_{x}^{2}\psi||_{X_{b}^{s}}||u_{k}||_{X_{b}^{s}}\nonumber
\\
\label{e416}& \leq & c\,A_{11}^{k}\,k!
\end{eqnarray}
\begin{eqnarray}
||t\,\partial_{x}^{3}\psi\cdot u_{k}||_{X_{b - 1}^{s - 1}} & \leq &
c\,||<D_{x,\,t}>^{3/2}t\,\partial_{x}^{3}\psi||_{X_{1 - b}^{|s| +
2|b - 1|}}||u_{k}||_{X_{b - 1}^{s}} \leq c\,||u_{k}||_{X_{b}^{s}}
\nonumber \\
\label{e417}& \leq & c\,A_{12}^{k}\,k!.
\end{eqnarray}
Hence, replacing \eqref{e413}, \eqref{e415},\eqref{e416} and
\eqref{e417} in \eqref{e414} we obtain that there exists a constant $c$
and $A_{14}$ such that
\begin{eqnarray}
\label{e418}||t\,\partial_{x}^{3}(\psi u_{k})||_{X_{b - 1}^{s - 1}}
 \leq  c\,A_{14}^{k}\cdot k! + c\,k!\sum_{k=k_{1}'' + k_{2}'' +
k_{3}''}\frac{1}{k_{1}''!}\,2^{k_{1}''}\,A_{9}^{k_{2}''}\cdot
A_{10}^{k_{3}''},\qquad k=0,\,1,\,2,\,\ldots
\end{eqnarray}
Therefore, replacing \eqref{e404}, \eqref{e405} and \eqref{e418} in
\eqref{e403} we obtain that there exists a constant $c$ and $A_{15}$
such that
\begin{eqnarray}
\lefteqn{||<D_{x,\,t}>^{3b}\psi
\,u_{k}||_{L_{t}^{2}(\mathbb{R}:\,H_{x}^{s - 1}(\mathbb{R}))} }
\nonumber \\
\label{e419}& \leq & c\,A_{15}^{k}\cdot k! + c\,k!\sum_{k=k_{1}'' +
k_{2}'' +
k_{3}''}\frac{1}{k_{1}''!}\,2^{k_{1}''}\,A_{9}^{k_{2}''}\cdot
A_{10}^{k_{3}''},\qquad k=0,\,1,\,2,\,\ldots
\end{eqnarray}
In a similar way, we obtain that there exists a constant $c$ and
$A_{16}$ such that
\begin{eqnarray}
\lefteqn{||<D_{x,\,t}>^{3b}\psi
\,v_{k}||_{L_{t}^{2}(\mathbb{R}:\,H_{x}^{s - 1}(\mathbb{R}))} }
\nonumber \\
\label{e420}& \leq & c\,A_{16}^{k}\cdot k! + c\,k!\sum_{k=k_{1}'' +
k_{2}'' +
k_{3}''}\frac{1}{k_{1}''!}\,2^{k_{1}''}\,A_{9}^{k_{2}''}\cdot
A_{10}^{k_{3}''},\qquad k=0,\,1,\,2,\,\ldots
\end{eqnarray}
Adding \eqref{e419} and \eqref{e420} we have
\begin{eqnarray}
\lefteqn{||<D_{x,\,t}>^{3b}\psi
\,u_{k}||_{L_{t}^{2}(\mathbb{R}:\,H_{x}^{s - 1}(\mathbb{R}))}
 + ||<D_{x,\,t}>^{3b}\psi \,v_{k}||_{L_{t}^{2}(\mathbb{R}:\,H_{x}^{s -
1}(\mathbb{R}))} }
\nonumber \\
& \leq & c\,A_{15}^{k}\cdot k! + c\,A_{16}^{k}\cdot k! +
c\,k!\sum_{k=k_{1}'' + k_{2}'' +
k_{3}''}\frac{1}{k_{1}''!}\,2^{k_{1}''}\;2\cdot A_{9}^{k_{2}''}\cdot
A_{10}^{k_{3}''} \nonumber \\
& \leq & c\,(A_{15}^{k} + A_{16}^{k})\cdot k! + c\,k!\sum_{k=k_{1}''
+ k_{2}'' + k_{3}''}\frac{1}{k_{1}''!}\,2^{k_{1}''}\;2\cdot
A_{9}^{k_{2}''}\cdot
A_{10}^{k_{3}''} \nonumber \\
\label{e421}& \leq & c\,A_{17}^{k}\cdot k! + c\,k!\sum_{k=k_{1}'' +
k_{2}'' + k_{3}''}\frac{1}{k_{1}''!}\,2^{k_{1}''}\;2\cdot
A_{9}^{k_{2}''}\cdot A_{10}^{k_{3}''}.
\end{eqnarray}
We estimate the last term on the right hand side of \eqref{e421}
\begin{eqnarray}
\sum_{k=k_{1}'' + k_{2}'' +
k_{3}''}\frac{1}{k_{1}''!}\,2^{k_{1}''}\;2\cdot A_{9}^{k_{2}''}\cdot
A_{10}^{k_{3}''} & = & \sum_{m=0}^{k}\sum_{j=0}^{m}\frac{1}{(m -
j)!}\;2^{(m - j)}\;2\cdot A_{9}^{j}\cdot
A_{10}^{k - m} \nonumber \\
& \leq & A_{10}^{k}\sum_{m=0}^{k}\sum_{j=0}^{m}\frac{1}{(m -
j)!}\;2\cdot \left(\frac{A_{9}}{2}\right)^{j}\cdot
\left(\frac{2}{A_{10}}\right)^{m} \nonumber \\
& \leq &
A_{10}^{k}\sum_{m=0}^{k}\sum_{j=0}^{m}\left[\left(\frac{A_{9}^{2}}{4}\right)^{j}
+ \left(\frac{4}{A_{10}^{2}}\right)^{m}\right] \nonumber \\
& \leq &
A_{10}^{k}\sum_{m=0}^{k}\sum_{j=0}^{m}j!\;\frac{\left(\frac{A_{9}^{2}}{4}\right)^{j}}{j!}
+
A_{10}^{k}\sum_{m=0}^{k}\sum_{j=0}^{m}m!\;\frac{\left(\frac{4}{A_{10}^{2}}\right)^{m}}{m!}
\nonumber \\
& \leq &
A_{10}^{k}\;k!\sum_{m=0}^{k}\sum_{j=0}^{m}\frac{\left(\frac{A_{9}^{2}}{4}\right)^{j}}{j!}
+
A_{10}^{k}\;k!\sum_{m=0}^{k}\sum_{j=0}^{m}\frac{\left(\frac{4}{A_{10}^{2}}\right)^{m}}{m!}
\nonumber \\
& \leq & e^{A_{9}^{2}/4}\;A_{10}^{k}\;k! +
e^{4/A_{10}^{2}}\;A_{10}^{k}\;k!
\nonumber \\
\label{e422}& \leq & c\;A_{10}^{k}\;k!.
\end{eqnarray}
Replacing \eqref{e422} in \eqref{e421} we obtain
\begin{eqnarray}
&  & ||<D_{x,\,t}>^{3b}\psi
\,u_{k}||_{L_{t}^{2}(\mathbb{R}:\,H_{x}^{s - 1}(\mathbb{R}))}
 + ||<D_{x,\,t}>^{3b}\psi \,v_{k}||_{L_{t}^{2}(\mathbb{R}:\,H_{x}^{s -
1}(\mathbb{R}))}\nonumber \\
&  & \leq c\,A_{17}^{k}\cdot k! + c\,A_{19}^{k}\cdot (k!)^{2}
\nonumber \\
&  & \leq c\,A_{17}^{k}\cdot (k!)^{2} + c\,A_{19}^{k}\cdot (k!)^{2}
\nonumber \\
\label{e423}&  & \leq c\,A_{20}^{k}\cdot (k!)^{2}
\end{eqnarray}
and the result follows. 
\begin{remark}
 a) For simplicity, we only illustrate the conclusion
for the case $s\geq -1/2 - \delta$ with $b=1/2 + \delta/3$ (for
small $\delta>0$) and the case $s=-3/4 + \delta$ and $b=7/12 -
\delta/3.$ If $s=-1/2 - \delta$ with $b=1/2 + \delta/3,$ the initial
data can involve Dirac's delta measure $\delta_{0}$ and the latter
is the
critical case of the local well-posedness.\\
b) The following inequality is simple to verify in both cases,
\begin{eqnarray*}
||\psi\,u_{k}||_{L_{x,\,t}^{2}(\mathbb{R}^{2})}\leq
||<D_{x}>^{3\,b}(\psi\,u_{k})||_{L_{t}^{2}(\mathbb{R}:\;H_{x}^{s -
1}(\mathbb{R}))}\leq
c\;||<D_{x,\,t}>^{3\,b}(\psi\,u_{k})||_{L_{t}^{2}(\mathbb{R}:\;H_{x}^{s
- 1}(\mathbb{R}))}.
\end{eqnarray*}
\end{remark}
\begin{proposition} {\it Under the same assumptions as in
Proposition 4.1, there exist positive constants $c$ and $A$ such
that}
\begin{eqnarray}
\label{e424}||\psi\,P^{k}u||_{H^{7/2}(\mathbb{R}^{2})}\leq
c\,A^{k}\,(k!)^{2},\quad k=0,\,1,\,2,\ldots\\
\label{e425}||\psi\,P^{k}v||_{H^{7/2}(\mathbb{R}^{2})}\leq
c\,A^{k}\,(k!)^{2},\quad k=0,\,1,\,2,\ldots
\end{eqnarray}
\end{proposition}
\noindent
{\it Proof.} We apply Lemma 2.4 to $\psi\,u_{k}\equiv \psi\,P^{k}u$
with $b=1$ and $r=0.$
\begin{eqnarray}
\lefteqn{||<D_{x,\,t}>^{3}\psi\,P^{k}u||_{L^{2}(\mathbb{R}:\,L_{x}^{2}(\mathbb{R}))}}
\nonumber
\\
\label{e426}&  & \leq
c\,\left(||\psi\,u_{k}||_{L^{2}(\mathbb{R}:\,L_{x}^{2}(\mathbb{R}))}
+
||t\partial_{x}^{3}(\psi\,u_{k})||_{L^{2}(\mathbb{R}:\,L_{x}^{2}(\mathbb{R}))}
+
||P^{3}(\psi\,u_{k})||_{L^{2}(\mathbb{R}:\,L_{x}^{2}(\mathbb{R}))}\right).
\end{eqnarray}
Therefore, if we wish to estimate the second term in the right hand
side of \eqref{e426} with the aid of the equation \eqref{e216}
\begin{eqnarray*}
\psi\,t\,\partial_{x}^{3}u_{k} = -\,\frac{1}{3}\,\psi\,Pu_{k} +
\frac{1}{3}\,\psi\,x\,\partial_{x}u_{k} + t\,\psi\,B_{k}
\end{eqnarray*}
it is necessary to estimate
$||\psi\,u_{k}||_{L_{t}^{2}(\mathbb{R}:\,H_{x}^{1}(\mathbb{R}))}$
which is not yet obtained. Hence, we start from the lower regularity
setting, i. e., applying \eqref{e227} in Lemma 2.4 to $\psi\,u_{k}$
with $\mu=1/2.$ Let $\psi_{1}$ be a smaller size of smooth cut-off
function with $\psi_{1}\leq\psi$ and $\psi_{1}=1$ around
$(x_{0},\,t_{0}).$ Applying \eqref{e227} a
$\psi\,u_{k}=\psi\,P^{k}u$ with $\mu=1/2$ we have
\begin{eqnarray}
\lefteqn{||<D_{x,\,t}>^{3}\psi_{1}\,P^{k}u||_{H^{-\,
5/2}(\mathbb{R}^{2})}\leq
c\,||<D_{x,\,t}>^{3}\psi_{1}\,P^{k}u||_{L^{2}(\mathbb{R}^{2})}}\nonumber\\
\label{e427}&  & \leq c\,\left(||\psi_{1}u_{k}||_{H^{-\,
5/2}(\mathbb{R}^{2})} + ||t\partial_{x}^{3}(\psi_{1}u_{k})||_{H^{-\,
5/2}(\mathbb{R}^{2})} + ||P^{3}(\psi_{1}u_{k})||_{H^{-\,
5/2}(\mathbb{R}^{2})}\right).
\end{eqnarray}
The first term on the right hand side of \eqref{e427} has already
been estimated. For the third term we have
\begin{eqnarray}
||P^{3}(\psi_{1}u_{k})||_{H^{-\,5/2}(\mathbb{R}^{2})} & \leq &
||P^{3}(\psi_{1}u_{k})||_{L_{x,\,t}^{2}(\mathbb{R}^{2})} \nonumber
\\
& = & \sum_{l=0}^{3}\frac{3!}{l!(3 - l)!}||(P^{3 -
l}\psi_{1})(P^{l}u_{k})||_{L_{x,\,t}^{2}(\mathbb{R}^{2})} \nonumber
\\
& \leq & \sum_{l=0}^{3}\frac{3!}{l!(3 - l)!}||P^{3 -
l}\psi_{1}||_{L_{x,\,t}^{\infty}(\mathbb{R}^{2})}
||P^{l}u_{k}||_{L_{x,\,t}^{2}(\mathbb{R}^{2})}
\nonumber
\\
& \leq & c\sum_{l=0}^{3}\frac{3!}{l!(3 - l)!}||P^{k +
l}u||_{L_{x,\,t}^{2}(\mathbb{R}^{2})} \nonumber
\\
& \leq & c\sum_{l=1}^{3}A_{1}^{k + l}k!
\ \leq \ c\,A_{2}^{k}k!  
\label{e428}\ \leq \ c\,A_{2}^{k}(k!)^{2}.
\end{eqnarray}
For the second term on the right side hand we use the same idea of
the remark above, using the dilation operator $P.$ Indeed,
\begin{eqnarray}
||t\,\partial_{x}^{3}(\psi_{1}\,u_{k})||_{H^{-5/2}}
 & \leq &
||\psi_{1}\,t\,\partial_{x}^{3}u_{k}||_{H^{-5/2}(\mathbb{R}^{2})} +
3\,||\partial_{x}^{2}(t\,\partial_{x}\psi_{1}\cdot
u_{k})||_{H^{-5/2}(\mathbb{R}^{2})}\nonumber \\
\label{e429}& + & 
3\;||\partial_{x}(t\,\partial_{x}^{2}\psi_{1}\cdot
u_{k})||_{H^{-5/2}(\mathbb{R}^{2})} +
||t\,(\partial_{x}^{3}\psi_{1})u_{k}||_{H^{-5/2}(\mathbb{R}^{2})}.
\end{eqnarray}
The last three term are bounded by the following:
\begin{eqnarray}
& & c\left(||\partial_{x}\psi_{1}||_{L_{x,\,t}^{\infty}(\mathbb{R}^{2})}
+ ||\partial_{x}^{2}\psi_{1}||_{L_{x,\,t}^{\infty}(\mathbb{R}^{2})}
+
||\partial_{x}^{3}\psi_{1}||_{L_{x,\,t}^{\infty}(\mathbb{R}^{2})}\right)
||\psi\,u_{k}||_{L_{x,\,t}^{2}(\mathbb{R}^{2})} \nonumber \\ & & \leq
c\,A_{3}^{k}\,k! 
\label{e430}\ \leq \ c\,A_{3}^{k}\,(k!)^{2}.
\end{eqnarray}
On the other hand, using
\begin{eqnarray}
||\psi_{1}\,t\,\partial_{x}^{3}u_{k}||_{H^{-5/2}(\mathbb{R}^{2})} &
\leq &
\frac{1}{3}\,||\psi_{1}\,Pu_{k}||_{L^{2}(\mathbb{R}:\,L_{x}^{2}(\mathbb{R}))}
+
\frac{1}{3}\,||x\,\psi_{1}\,\partial_{x}u_{k}||_{H^{-5/2}(\mathbb{R}^{2})}
+ ||t\,\psi_{1}\,B_{k}||_{H^{-5/2}(\mathbb{R}^{2})}\nonumber
\\
\label{e431}& = & F_{1} + F_{2} + F_{3}.
\end{eqnarray}
Thus
\begin{eqnarray}
F_{1} & \leq &
c\,||\psi_{1}||_{L_{x,\,t}^{\infty}(\mathbb{R}^{2})}\,
||\psi\,P^{k + 1}u||_{L_{x,\,t}^{2}(\mathbb{R}^{2})}
\ \leq \ c\,||\psi\,P^{k +
1}v||_{L_{x,\,t}^{2}(\mathbb{R}^{2})}\nonumber \\
& \leq & c\,A_{4}^{k + 1}(k + 1)!
\ \leq \ c\,A_{5}^{k}\,k! 
\label{e432}\ \leq \ c\,A_{5}^{k}\,(k!)^{2},
\end{eqnarray}
\begin{eqnarray}
F_{2} & \leq &
||x\,\psi_{1}\,\partial_{x}v_{k}||_{L^{2}(\mathbb{R}:\,H_{x}^{-1}(\mathbb{R}))}\nonumber
\\
& \leq &
||\partial_{x}(x\,\psi_{1}\,v_{k})||_{L^{2}(\mathbb{R}:\,H_{x}^{-1}(\mathbb{R}))}
+
||\partial_{x}(x\,\psi_{1})\psi\,v_{k}||_{L^{2}(\mathbb{R}:\,H_{x}^{-1}(\mathbb{R}))}
\nonumber \\
& \leq & ||x\,\psi_{1}\,v_{k}||_{L_{x,\,t}^{2}(\mathbb{R}^{2})} +
||\partial_{x}(x\,\psi_{1})||_{L_{x,\,t}^{\infty}(\mathbb{R}^{2})}
||\psi\,v_{k}||_{L_{x,\,t}^{2}(\mathbb{R}^{2})}\nonumber \\
& \leq & \left(||x\,\psi_{1}||_{L_{x,\,t}^{\infty}(\mathbb{R}^{2})}
+
||\partial_{x}(x\,\psi_{1}))||_{L_{x,\,t}^{\infty}(\mathbb{R}^{2})}\right)
||\psi\,v_{k}||_{L_{x,\,t}^{2}(\mathbb{R}^{2})}\nonumber \\
& \leq & c\,A_{6}^{k}\,k! 
\label{e433}\ \leq \ c\,A_{6}^{k}\,(k!)^{2}.
\end{eqnarray}
Using Lemma 2.5(case $\sigma=-5/2,$ $s=5,$ $r=-5/2$) 
\begin{eqnarray*}
F_{3}  = 
||t\,\psi_{1}\,B_{k}||_{H^{-5/2}(\mathbb{R}^{2})}
 \leq  c_{1}\,||\psi_{1}||_{H^{5}(\mathbb{R}^{2})}
||\psi^{2}\,B_{x}||_{H^{-5/2}(\mathbb{R}^{2})}
\end{eqnarray*}
and replacing $B_{k}$ by \eqref{e214}, we have
\begin{eqnarray}
F_{3} & \leq & c_{1}\,\frac{|a|}{2}\sum_{k=k_{1} + k_{2} +
k_{3}}\frac{k!}{k_{1}!\,k_{2}!\,k_{3}!}2^{k_{1}}||\psi\,u_{k_{2}}
\psi\,u_{k_{3}}||_{H^{-3/2}(\mathbb{R}^{2})}\nonumber
\\
&  & +\;c_{1}\,\frac{|b|}{2}\sum_{k=k_{1}' + k_{2}' +
k_{3}'}\frac{k!}{k_{1}'!\,k_{2}'!\,k_{3}'!}2^{k_{1}'}||\psi\,v_{k_{2}'}
\psi\,v_{k_{3}'}||_{H^{-3/2}(\mathbb{R}^{2})}\nonumber
\\
&  & +\;c_{1}\,|c|\sum_{k=k_{1}'' + k_{2}'' +
k_{3}''}\frac{k!}{k_{1}''!\,k_{2}''!\,k_{3}''!}2^{k_{1}''}||\psi\,u_{k_{2}''}
\psi\,v_{k_{3}''}||_{H^{-3/2}(\mathbb{R}^{2})}\nonumber
\\
& \leq & c_{1}\,\frac{|a|}{2}\sum_{k=k_{1} + k_{2} +
k_{3}}\frac{k!}{k_{1}!\,k_{2}!\,k_{3}!}2^{k_{1}}
||\psi\,u_{k_{2}}||_{L^{2}(\mathbb{R}^{2})}
||\psi\,u_{k_{3}}||_{L^{2}(\mathbb{R}^{2})}\nonumber
\\
&  & +\;c_{1}\,\frac{|b|}{2}\sum_{k=k_{1}' + k_{2}' +
k_{3}'}\frac{k!}{k_{1}'!\,k_{2}'!\,k_{3}'!}2^{k_{1}'}
||\psi\,v_{k_{2}'}||_{L^{2}(\mathbb{R}^{2})}
||\psi\,v_{k_{3}'}||_{L^{2}(\mathbb{R}^{2})}\nonumber
\\
&  & +\;c_{1}\,|c|\sum_{k=k_{1}'' + k_{2}'' +
k_{3}''}\frac{k!}{k_{1}''!\,k_{2}''!\,k_{3}''!}2^{k_{1}''}
||\psi\,u_{k_{2}''}||_{L^{2}(\mathbb{R}^{2})}
||\psi\,v_{k_{3}''}||_{L^{2}(\mathbb{R}^{2})}\nonumber
\\
& \leq & c_{1}\,\frac{|a|}{2}\sum_{k=k_{1} + k_{2} +
k_{3}}\frac{k!}{k_{1}!\,k_{2}!\,k_{3}!}2^{k_{1}}\,
A_{7}^{k_{2}}k_{2}!A_{7}^{k_{3}}k_{3}!\nonumber
\\
&  & +\;c_{1}\,\frac{|b|}{2}\sum_{k=k_{1}' + k_{2}' +
k_{3}'}\frac{k!}{k_{1}'!\,k_{2}'!\,k_{3}'!}2^{k_{1}'}\,
A_{8}^{k_{2}'}k_{2}'!A_{8}^{k_{3}'}k_{3}'!\nonumber
\\
&  & +\;c_{1}\,|c|\sum_{k=k_{1}'' + k_{2}'' +
k_{3}''}\frac{k!}{k_{1}''!\,k_{2}''!\,k_{3}''!}2^{k_{1}''}\,
A_{9}^{k_{2}''}k_{2}''!A_{10}^{k_{3}''}k_{3}''!\nonumber
\\
& \leq & c_{1}\,\frac{|a|}{2}\,k!\sum_{k=k_{1} + k_{2} +
k_{3}}\frac{2^{k_{1}}}{k_{1}!}\, A_{7}^{k_{2} + k_{3}} +
c_{1}\,\frac{|b|}{2}\,k!\sum_{k=k_{1}' + k_{2}' +
k_{3}'}\frac{2^{k_{1}'}}{k_{1}'!}\, A_{8}^{k_{2}' + k_{3}'}\nonumber
\\
&  & +\;c_{1}\,|c|\,k!\sum_{k=k_{1}'' + k_{2}'' +
k_{3}''}\frac{2^{k_{1}''}}{k_{1}''!}\,
A_{9}^{k_{2}''}A_{10}^{k_{3}''},\nonumber
\end{eqnarray}
and then
\begin{eqnarray}
F_3 & \leq &
c_{1}\,\frac{|a|}{2}\,k!\,A_{7}^{k}\sum_{k_{1}=0}^{k}\sum_{k_{2}=0}^{k
- k_{1}}\frac{2^{k_{1}}}{k_{1}!}\, A_{7}^{-\,k_{1}} +
c_{1}\,\frac{|b|}{2}\,k!\,A_{8}^{k}\sum_{k_{1}=0}^{k}\sum_{k_{2}'=0}^{k
-
k_{1}'}\frac{2^{k_{1}'}}{k_{1}'!}\, A_{8}^{-\,k_{1}'}\nonumber \\
&  & +\;c_{1}\,|c|\,k!\sum_{k=k_{1}'' + k_{2}'' +
k_{3}''}\frac{2^{k_{1}''}}{k_{1}''!}\,
A_{9}^{k_{2}''}A_{10}^{k_{3}''}\nonumber
\\
& \leq &
c_{1}\,\frac{|a|}{2}\,e^{2/A_{7}}\,A_{7}^{k}\,(k + 1)! +
c_{1}\,\frac{|b|}{2}\,e^{3/A_{8}}\,A_{8}^{k}\,(k + 1)!
\nonumber\\ & &
+
\;c_{1}\,|c|\,k!\sum_{k=k_{1}'' + k_{2}'' +
k_{3}''}\frac{2^{k_{1}''}}{k_{1}''!}\,
A_{9}^{k_{2}''}A_{10}^{k_{3}''}.
\label{e434}
\end{eqnarray}
Replacing \eqref{e430}, \eqref{e435} and \eqref{e429} in
\eqref{e431} we obtain
\begin{eqnarray}
\lefteqn{||\psi_{1}\,t\,\partial_{x}^{3}u_{k}||_{H^{-5/2}(\mathbb{R}^{2})}
} \nonumber \\
\label{e435}& \leq & c_{2}\,A_{11}^{k}\,k! +
c_{1}\,|c|\,k!\sum_{k=k_{1}'' + k_{2}'' +
k_{3}''}\frac{2^{k_{1}''}}{k_{1}''!}\,
A_{9}^{k_{2}''}A_{10}^{k_{3}''},\quad k=0,\,1,\,2,\,\ldots
\end{eqnarray}
Replacing \eqref{e430} and \eqref{e435} in \eqref{e429}
\begin{eqnarray}
\lefteqn{||t\,\partial_{x}^{3}(\psi_{1}\,u_{k})||_{H^{-5/2}(\mathbb{R}^{2})}
}\nonumber \\
\label{e436}& \leq & c_{3}\,A_{12}^{k}\,k! +
c_{1}\,|c|\,k!\sum_{k=k_{1}'' + k_{2}'' +
k_{3}''}\frac{2^{k_{1}''}}{k_{1}''!}\,
A_{9}^{k_{2}''}A_{10}^{k_{3}''},\quad k=0,\,1,\,2,\,\ldots
\end{eqnarray}
Now replacing \eqref{e428} and \eqref{e436} in \eqref{e427} we
obtain
\begin{eqnarray}
\lefteqn{||<D_{x,\,t}>^{3}\psi\,u_{k}||_{H^{-\,
5/2}(\mathbb{R}^{2})}} \nonumber \\
\label{e437}& \leq & c_{4}\,A_{13}^{k}\,k! +
c_{1}\,|c|\,k!\sum_{k=k_{1}'' + k_{2}'' +
k_{3}''}\frac{2^{k_{1}''}}{k_{1}''!}\,
A_{9}^{k_{2}''}A_{10}^{k_{3}''},\quad k=0,\,1,\,2,\,\ldots
\end{eqnarray}
In particular
\begin{eqnarray}
\lefteqn{||\psi\,u_{k}||_{H^{1/2}(\mathbb{R}^{2})}} \nonumber \\
\label{e438}& \leq & c_{5}\,A_{14}^{k}\,k! +
c_{1}\,|c|\,k!\sum_{k=k_{1}'' + k_{2}'' +
k_{3}''}\frac{2^{k_{1}''}}{k_{1}''!}\,
A_{9}^{k_{2}''}A_{10}^{k_{3}''},\quad k=0,\,1,\,2,\,\ldots
\end{eqnarray}
Using a similar argument as above for
$||<D_{x,\,t}>^{3}\psi\,P^{k}u||_{H^{-\, 3/2}(\mathbb{R}^{2})}$ with
$\mu = 3/2$ in \eqref{e227} and replacing the support of the cut-off
function $\psi_{\epsilon}$ we obtain
\begin{eqnarray}
\lefteqn{||\psi\,u_{k}||_{H^{3/2}(\mathbb{R}^{2})} }\nonumber \\
\label{e439} & \leq &  c_{5}\,A_{14}^{k}\,k! +
c_{1}\,|c|\,k!\sum_{k=k_{1}'' + k_{2}'' +
k_{3}''}\frac{2^{k_{1}''}}{k_{1}''!}\,
A_{9}^{k_{2}''}A_{10}^{k_{3}''},\quad k=0,\,1,\,2,\,\ldots
\end{eqnarray}
In a similar way we have
\begin{eqnarray}
\lefteqn{||\psi\,v_{k}||_{H^{3/2}(\mathbb{R}^{2})} }\nonumber \\
\label{e440} & \leq &  c_{5}\,A_{15}^{k}\,k! +
c_{1}\,|\widetilde{c}|\,k!\sum_{k=k_{1}'' + k_{2}'' +
k_{3}''}\frac{2^{k_{1}''}}{k_{1}''!}\,
A_{9}^{k_{2}''}A_{10}^{k_{3}''},\quad k=0,\,1,\,2,\,\ldots
\end{eqnarray}
Adding \eqref{e439} with \eqref{e440} and performing straightforward
calculations as \eqref{e422} we obtain
\begin{eqnarray}
\label{e441}||\psi\,u_{k}||_{H^{3/2}(\mathbb{R}^{2})} +
||\psi\,v_{k}||_{H^{3/2}(\mathbb{R}^{2})} \leq
C\,A^{k}\,(k!)^{2},\quad k=0,\,1,\,2,\,\ldots
\end{eqnarray}
To obtain the estimate for $||\psi\,P^{k}u||_{H^{
7/2}(\mathbb{R}^{2})}$ and $||\psi\,P^{k}v||_{H^{
7/2}(\mathbb{R}^{2})}$ we repeat the above method with $\mu=7/2.$
\begin{proposition} 
{\it Suppose that}
\begin{eqnarray}
\label{e442}||\psi\,u_{k}||_{H^{7/2}(\mathbb{R}^{2})} \leq
c\,A_{1}^{k}\,(k!)^{2},\quad k=0,\,1,\,2,\,\ldots \\
\label{e443}||\psi\,v_{k}||_{H^{7/2}(\mathbb{R}^{2})} \leq
c\,A_{2}^{k}\,(k!)^{2},\quad k=0,\,1,\,2,\,\ldots
\end{eqnarray}
{\it then we have}
\begin{eqnarray}
\sup_{t\in [t_{0} - \epsilon,\,t_{0} + \epsilon]}
\label{e444}||(t^{1/3}\partial_{x})P^{k}u||_{H^{1}(x_{0} -
\epsilon,\,x_{0} + \epsilon)} \leq c_{1}\,A_{3}^{k + l}\,[\,(k +
l)!]^{2},\quad k,\,l=0,\,1,\,2,\,\ldots \\
\sup_{t\in [t_{0} - \epsilon,\,t_{0} + \epsilon]}
\label{e445}||(t^{1/3}\partial_{x})P^{k}v||_{H^{1}(x_{0} -
\epsilon,\,x_{0} + \epsilon)} \leq c_{1}\,A_{4}^{k + l}\,[\,(k +
l)!\,]^{2},\quad k,\,l=0,\,1,\,2,\,\ldots
\end{eqnarray}
{\it where $\epsilon>0$ is so small that $\psi\equiv 1$ near
$I=(x_{0} - \epsilon,\,x_{0} + \epsilon)\times (t_{0} - \epsilon,\,t_{0} + \epsilon).$}\\
\end{proposition}
\noindent
{\it Proof.} Let $I_{t_{0}}=(t_{0} - \epsilon,\,t_{0} + \epsilon)$
and $I_{x_{0}}=(x_{0} - \epsilon,\,x_{0} + \epsilon),$ then we have
$I=I_{x_{0}}\times I_{t_{0}}.$ For any fixed $t\in I_{x_{0}},$ let
${\cal L}=t^{1/3}\partial_{x}.$ We show that for some positive
constants $c$ and $A_{0}$ the following inequality holds
\begin{eqnarray}
\label{e446}||{\cal L}^{l}P^{k}u||_{H_{x}^{1}(I_{x_{0}})}\leq
c\,A_{0}^{k + l}[\,(k + l)!\,]^{2}, \quad \forall
\,k,\;\forall\,l=0,\,1,\,2,\,\ldots
\end{eqnarray}
Now, let use induction over $l.$ By the trace theorem, we have
\begin{eqnarray}
||{\cal L}^{l}P^{k}u||_{H_{x}^{1}(I_{x_{0}})} & \leq &
||t^{l/3}\,\partial_{x}^{l}P^{k}u(t)||_{H_{x}^{1}(I_{x_{0}})}
\ \leq \ (t_{0} +
\epsilon)^{l/3}||\partial_{x}^{l}P^{k}u||_{H^{3/2}(I_{x_{0}}\times
I_{t_{0}})}
\nonumber \\
& \leq & (t_{0} + \epsilon)^{l/3}||P^{k}u||_{H^{7/2}(I_{x_{0}}\times
I_{t_{0}})}
\ \leq \ (t_{0} +
\epsilon)^{l/3}||\psi\,P^{k}u||_{H^{7/2}(\mathbb{R}^{2})}
\nonumber \\
& \leq & (t_{0} +
\epsilon)^{l/3}\,c_{1}\,A_{1}^{k}\,k!
\ \leq \ (t_{0} +
\epsilon)^{l/3}\,c_{1}\,A_{0}^{k + l}\,(k + l) \nonumber \\
\label{e447}& \leq & (t_{0} + \epsilon)^{l/3}\,c_{1}\,A_{0}^{k +
l}\,[\,(k + l)!\,]^{2}.
\end{eqnarray}
where we take $c=(t_{0} + \epsilon)^{l/3}c_{1}$ and
$A_{0}=\mbox{max}\{1,\,A_{1}\}.$ Hence, in the case $l=0,\,1,\,2,$
it is easy to show that \eqref{e446} follows directly from the
assumption. \\
Now, we assume that \eqref{e446} is true to $l\geq 2.$ Applying
$P^{k}$ to the  equation \eqref{e207}, we have
\begin{eqnarray*}
\partial_{t}(P^{k}u) + \partial_{x}^{3}(P^{k}u)
& = & L\,P^{k}u  \\
& = & (P + 3)^{k}Lu  \\
& = & (P + 3)^{k}(\partial_{t}u + \partial_{x}^{3}u) \\
& = & -(P + 3)^{k}\left[\frac{a}{2}\;\partial_{x}(u^{2}) +
\frac{b}{2}\;\partial_{x}(v^{2}) + c\,\partial_{x}(u\,v)\right]
 \\
& = & -\,\frac{a}{2}\,(P + 3)^{k}\partial_{x}(u^{2}) -
\frac{b}{2}\,(P +
3)^{k}\partial_{x}(v^{2}) - c\,(P + 3)^{k}\partial_{x}(u\,v) \\
& = & -\,\frac{a}{2}\;\partial_{x}(P + 2)^{k}(u^{2}) -
\frac{b}{2}\;\partial_{x}(P + 2)^{k}(v^{2}) - c\;\partial_{x}(P +
2)^{k}(u\,v)
\end{eqnarray*}
such that
\begin{eqnarray}
\label{e448}t\,\partial_{t}(P^{k}u) + t\,\partial_{x}^{3}(P^{k}u) =
-\,\frac{a}{2}\;t\,\partial_{x}(P + 2)^{k}(u^{2}) -
\frac{b}{2}\;t\,\partial_{x}(P + 2)^{k}(v^{2}) -
c\;t\,\partial_{x}(P + 2)^{k}(u\,v).
\end{eqnarray}
Moreover, $P=3\,t\,\partial_{t} + x\,\partial_{x}.$ Then
\begin{eqnarray}
\label{e449}t\,\partial_{t}(P^{k}u) = \frac{1}{3}\,P^{k + 1}u -
\frac{1}{3}\,x\,\partial_{x}(P^{k}u).
\end{eqnarray}
Replacing \eqref{e449} in  \eqref{e448} we obtain
\begin{eqnarray}
{\cal L}^{3}P^{k}u & = & t\,\partial_{x}^{3}(P^{k}u)
 = -\;\frac{1}{3}\,P^{k + 1}u + \frac{1}{3}\,x\,\partial_{x}(P^{k}u) \nonumber \\
\label{e450}&  & -\,\frac{a}{2}\;t\,\partial_{x}(P + 2)^{k}(u^{2}) -
\frac{b}{2}\;t\,\partial_{x}(P + 2)^{k}(v^{2}) -
c\;t\,\partial_{x}(P + 2)^{k}(u\,v).
\end{eqnarray}
Hence, applying ${\cal L}^{l - 2}$ we have
\begin{eqnarray}
||{\cal L}^{l + 1}P^{k}u||_{H_{x}^{1}(I_{x_{0}})} & = &
||{\cal L}^{l - 2}{\cal L}^{3}P^{k}u||_{H_{x}^{1}(I_{x_{0}})}\nonumber \\
& \leq & \frac{1}{3}\,||{\cal L}^{l - 2}\,P^{k +
1}u||_{H_{x}^{1}(I_{x_{0}})} + \frac{1}{3}\,||{\cal L}^{l -
2}\,x\,\partial_{x}(P^{k}u)||_{H_{x}^{1}(I_{x_{0}})}
\nonumber \\
&  & +\;\frac{|a|}{2}\,||t\,{\cal L}^{l - 2}\,\partial_{x}(P +
2)^{k}(u^{2})||_{H_{x}^{1}(I_{x_{0}})} + \frac{|b|}{2}\,||t\,{\cal
L}^{l - 2}\,\partial_{x}(P + 2)^{k}(v^{2})||_{H_{x}^{1}(I_{x_{0}})}
\nonumber \\
&  & \;+\,|c|\,||t\,{\cal L}^{l - 2}\,\partial_{x}(P +
2)^{k}(u\,v)||_{H_{x}^{1}(I_{x_{0}})} \nonumber \\
\label{e451}& = & F_{1} + F_{2} + F_{3} + F_{4} + F_{5}.
\end{eqnarray}
Using the induction assumption, we obtain
\begin{eqnarray}
\label{e452}F_{1} \leq \frac{1}{3}\,c_{1}\,A_{14}^{k + l + 1}(k + l
+ 1)!.
\end{eqnarray}
We estimate the term ${\cal L}^{l - 2}(x\,\partial_{x})$ for $l\geq
3.$ Let $r=l - 2,$ then we estimate ${\cal L}^{r}(x\,\partial_{x})$
for $r\geq 1.$
\begin{eqnarray}
\label{e453}\partial_{x}^{r}(x\,\partial_{x}) =
\sum_{k=0}^{r}{r\choose k}\,\partial_{x}^{r - k}(\,x\,)\cdot
\partial_{x}^{k}(\,\partial_{x}\,).
\end{eqnarray}
But
\[ \partial_{x}^{r - k}(\,x\,) = \left \{ \begin{array}{ll}
1 &
\mbox{if}\quad k=r - 1\\
0 & \mbox{if}\quad k\leq r - 2
\end{array}
\right. \]
then in \eqref{e453} we obtain
\begin{eqnarray*}
\partial_{x}^{r}(x\,\partial_{x}) & = & r\,\partial_{x}^{r - 1}(\,\partial_{x}\,) +
x\,\partial_{x}^{r}(\,\partial_{x}\,)
\ = \ r\,\partial_{x}^{r} +
x\,\partial_{x}(\,\partial_{x}^{r}\,)\\
& = & (l - 2)\,\partial_{x}^{(l - 2)} +
x\,\partial_{x}(\,\partial_{x}^{(l - 2)}\,),
\end{eqnarray*}
that is,
$\displaystyle {\cal L}^{l - 2}(x\,\partial_{x}) = x\,\partial_{x}{\cal
L}^{l - 2} + (l - 2)\,{\cal L}^{l - 2},$ for  $l\geq 3.
$
For $F_{2}$ we have
\begin{eqnarray}
F_{2} & \leq & ||x\,\partial_{x}{\cal L}^{l -
2}P^{k}u||_{H_{x}^{1}(I_{x_{0}})} + (l - 2)\,||{\cal L}^{l -
2}P^{k}u||_{H_{x}^{1}(I_{x_{0}})} \nonumber \\
& \leq & ||x\,t^{-1/3}{\cal L}^{l -
1}P^{k}u||_{H_{x}^{1}(I_{x_{0}})} + (l - 2)\,||{\cal L}^{l -
2}P^{k}u||_{H_{x}^{1}(I_{x_{0}})} \nonumber \\
& \leq & c\,(t_{0} - \epsilon)\,(|x_{0}| + \epsilon + 1)\,||{\cal
L}^{l - 1}P^{k}u||_{H_{x}^{1}(I_{x_{0}})} + (l - 2)\,||{\cal L}^{l -
2}P^{k}u||_{H_{x}^{1}(I_{x_{0}})} \nonumber \\
& \leq & (t_{0} - \epsilon)^{-1/3}(|x_{0} + \epsilon +
1)\,c_{1}\,A_{14}^{k + l - 1}(k + l - 1)! + c_{1}\,A_{14}^{k + l -
1}(l - 2)\,(k + l - 1)!\nonumber \\
\label{e455}& \leq & \frac{1}{3}\,c_{1}\,A_{14}^{k + l + 1}(k + l +
1)!
\end{eqnarray}
where we take $A_{14}$ larger than $(t_{0} -
\epsilon)^{-1/3}(|x_{0}| + \epsilon + 1)$ and 3. Using that (${\cal
L}=t^{1/3}\,\partial_{x}^{3}$)
\begin{eqnarray*}
t\,{\cal L}^{l - 2}\partial_{x} = t\,t^{(l - 2)/3}\partial_{x}^{(l -
2)}\partial_{x} = t\,t^{-1/3}\,t^{(l - 1)/3}\partial_{x}^{(l - 1)} =
t^{2/3}\,{\cal L}^{l - 1},
\end{eqnarray*}
we have
\begin{eqnarray*}
F_{3} & = & \frac{|a|}{2}\,
||t^{2/3}\,{\cal L}^{l - 1}\,(P + 2)^{k}(u^{2})||_{H_{x}^{1}(I_{x_{0}})}\nonumber \\
& \leq & \frac{|a|}{2}\,(t_{0} + \epsilon)^{2/3}\sum_{l - 1=l_{1} +
l_{2}}\sum_{k=k_{1} + k_{2} + k_{3}}\frac{(l - 1)!}{l_{1}!l_{2}!}\,
\frac{k!}{k_{1}!k_{2}!k_{3}!}\,2^{k_{3}} \\
&  & \times \;c_{2}\,||{\cal
L}^{l_{1}}P^{k_{1}}u||_{H_{x}^{1}(I_{x_{0}})}\,||{\cal
L}^{l_{2}}P^{k_{2}}u||_{H_{x}^{1}(I_{x_{0}})}.
\end{eqnarray*}
Using the induction assumption
\begin{eqnarray*}
F_{3} & \leq & \frac{|a|}{2}\,(t_{0} + \epsilon)^{2/3}\sum_{l -
1=l_{1} + l_{2}} \sum_{k=k_{1} + k_{2} +
k_{3}}c_{2}\,c_{1}^{3}\,k!\,(l - 1)!\,
\frac{2^{k_{3}}}{k_{3}!} \\
&  & \times \;\frac{(l_{1} + k_{1})!}{l_{1}!\,k_{1}!}\, \frac{(l_{2}
+ k_{2})!}{l_{2}!\,k_{2}!}\,
A_{14}^{k + l - 1} \\
& \leq & \frac{|a|}{2}\,(t_{0} +
\epsilon)^{2/3}\,c_{2}\,c_{1}^{3}\,(l + k - 1)!\, A_{14}^{k + l - 1}
\sum_{l - 1=l_{1} + l_{2}}\sum_{k=k_{1} + k_{2} + k_{3}}
\frac{2^{k_{3}}}{k_{3}!} \\
&  & \times \;\frac{(l_{1} + k_{1})!}{l_{1}!\,k_{1}!}\, \frac{(l_{2}
+ k_{2})!}{l_{2}!\,k_{2}!}\,\frac{k!\,(l - 1)!}{(l + k - 1)!}.
\end{eqnarray*}
Using that
\begin{eqnarray*}
\sum_{l - 1=l_{1} + l_{2}}\sum_{k=k_{1} + k_{2} + k_{3}}
\frac{2^{k_{3}}}{k_{3}!} \;\frac{(l_{1} + k_{1})!}{l_{1}!\,k_{1}!}\,
\frac{(l_{2} + k_{2})!}{l_{2}!\,k_{2}!}\,\frac{k!\,(l - 1)!}{(l + k
- 1)!}\leq e^{2}(l + k)!
\end{eqnarray*}
we obtain
\begin{eqnarray}
F_{3}  \ \leq  \ (t_{0} + \epsilon)^{2/3}\,c_{2}\,c_{1}^{3}\,e^{2}\,(l
+ k)!\,
A_{14}^{k + l - 1}
\label{e456}\  \leq  \ \frac{1}{3}\,c_{1}\,A_{14}^{k + l + 1}\,(k + l
+ 1)!
\end{eqnarray}
where we take $A_{14}$ larger than $(t_{0} -
\epsilon)^{-1/3}\,c_{2}\,c_{1}^{2}\,e^{2},$ and 3. In a similar way
\begin{eqnarray}
\label{e457} F_{4} \leq \frac{1}{3}\,c_{3}\,A_{15}^{k + l + 1}\,(k +
l + 1)!
\end{eqnarray}
where we take $A_{15}$ larger than $(t_{0} -
\epsilon)^{-1/3}\,c_{4}\,c_{3}^{2}\,e^{2},$ and 3. Finally, in a
similar way
\begin{eqnarray}
\label{e458} F_{5} \leq \frac{1}{3}\,c_{6}\,A_{16}^{k + l + 1}\,(k +
l + 1)!
\end{eqnarray}
where we take $A_{16}$ larger than $(t_{0} -
\epsilon)^{-1/3}\,c_{6}\,c_{5}^{2}\,e^{2},$ and 3. Therefore, from
\eqref{e452}, \eqref{e455}, \eqref{e456}, \eqref{e457} and
\eqref{e458} we obtain
\begin{eqnarray}
\label{e459}||{\cal L}^{l + 1}P^{k}u||_{H_{x}^{1}(I_{x_{0}})} \leq
c_{7}\,A_{17}^{k + l + 1}\,(k + l + 1)!.
\end{eqnarray}
In a similar way, we obtain
\begin{eqnarray}
\label{e460}||{\cal L}^{l + 1}P^{k}v||_{H_{x}^{1}(I_{x_{0}})} \leq
c_{7}\,A_{17}^{k + l + 1}\,(k + l + 1)!,
\end{eqnarray}
and the result follows. 
\begin{proposition} {\it Suppose that there exists a positive
constants $c_{1},\,c_{2}$ and $A_{14},\,A_{15}$ such that}
\begin{eqnarray}
\label{e461}\sup_{t\in [t_{0} - \epsilon,\,t_{0} + \epsilon]}
||\partial_{x}^{l}P^{k}u||_{H_{x}^{1}(x_{0} - \epsilon,\,x_{0} +
\epsilon)} \leq c_{1}\,A_{14}^{k + l}\,[\,(k + l)!\,]^{2},\quad
k,\,l=0,\,1,\,2,\,\ldots \\
\label{e462}\sup_{t\in [t_{0} - \epsilon,\,t_{0} + \epsilon]}
||\partial_{x}^{l}P^{k}v||_{H_{x}^{1}(x_{0} - \epsilon,\,x_{0} +
\epsilon)} \leq c_{2}\,A_{15}^{k + l}\,[\,(k + l)!\,]^{2},\quad
k,\,l=0,\,1,\,2,\,\ldots
\end{eqnarray}
{\it Then we have respectively}
\begin{eqnarray}
\label{e463}\sup_{t\in [t_{0} - \epsilon,\,t_{0} + \epsilon]}
||\partial_{t}^{m}\partial_{x}^{l}u||_{H_{x}^{1}(x_{0} -
\epsilon,\,x_{0} + \epsilon)} \leq c_{3}\,A_{16}^{m + l}\,[\,(m +
l)!\,]^{2},\quad m,\,l=0,\,1,\,2,\,\ldots \\
\label{e464}\sup_{t\in [t_{0} - \epsilon,\,t_{0} + \epsilon]}
||\partial_{t}^{m}\partial_{x}^{l}v||_{H_{x}^{1}(x_{0} -
\epsilon,\,x_{0} + \epsilon)} \leq c_{4}\,A_{17}^{m + l}\,[\,(m +
l)!\,]^{2},\quad m,\,l=0,\,1,\,2,\,\ldots
\end{eqnarray}
{\it where $c_{3},\,c_{4}$ and $A_{16},\,A_{17}$ only depend on
$c_{1},\,c_{2}$ and $A_{14},\,A_{15},$ respectively and $\epsilon,$
$(x_{0},\,t_{0}).$}
\end{proposition}
\noindent
{\it Proof.} Using the idea of Proposition 4.3, we fix $t\in
I_{x_{0}}.$ First we show that for some positive constants $c_{3},$
$A_{16}$ and $B_{16}$
\begin{eqnarray}
\label{e465}||
(x\,\partial_{x})^{m}\,\partial_{x}^{l}P^{k}v||_{H_{x}^{1}(I_{x_{0}})}
\leq c_{3}\,A_{16}^{k + m + l}\,B_{16}^{m}(k + m + l)!,\quad
k,\,m,\,l=0,\,1,\,2,\,\ldots
\end{eqnarray}
We use induction. Suppose that \eqref{e465} is true for $m.$\\
\begin{eqnarray}
\lefteqn{||(x\,\partial_{x})^{m +
1}\,\partial_{x}^{l}P^{k}v||_{H_{x}^{1}(I_{x_{0}})}}
\nonumber \\
& = &
||(x\,\partial_{x})\,(x\,\partial_{x})^{m}\,\partial_{x}^{l}P^{k}v||_{H_{x}^{1}(I_{x_{0}})}\nonumber \\
& \leq & (|x_{0}| + \epsilon + 1)\,||(x\,\partial_{x} + I)^{m}\,
\partial_{x}^{l + 1}P^{k}v||_{H_{x}^{1}(I_{x_{0}})}\nonumber \\
& \leq & c(|x_{0}|,\,\epsilon)\sum_{j=1}^{m}{m\choose j}
||(x\,\partial_{x})^{j}\,\partial_{x}^{l + 1}P^{k}v||_{H_{x}^{1}(I_{x_{0}})}\nonumber \\
& \leq & c\sum_{j=1}^{m}{m\choose j}c_{3}\,
A_{16}^{k + l + j + 1}\,B_{16}^{j}(k + l + j + 1)!\nonumber \\
& \leq & c_{3}\, A_{16}^{k + l + m + 1}\,B_{16}^{m}(k + l + m + 1)!
\sum_{j=1}^{m}\frac{(A_{16}\,B_{16})^{-(m - j)}}{(m - j)!}\;
\frac{m!}{j!}\;\frac{(k + l + j + 1)!}{(k + l + m + 1)!}\nonumber \\
\label{e466}& \leq & e^{-\,A_{16}\,B_{16}}c_{3}\,A_{16}^{k + l + m +
1}\,B_{16}^{m}(k + l + m + 1)!
\end{eqnarray}
where we take $B_{16}$ so large that $B_{16}\geq \mbox{max}\{|x_{0}|
+ \epsilon + 1,\,1\}.$  We show that for some positive constants
$c_{4},$ $A_{17}$ we have
\begin{eqnarray*}
||(t\,\partial_{t})^{m}\,\partial_{x}^{l}u||_{H_{x}^{1}(I_{x_{0}})}\leq
c_{4}\,A_{17}^{l + m}\,(l + m)!,\quad l,\,m=0,\,1,\,2,\,\ldots
\end{eqnarray*}
Using that $t\,\partial_{t}=\frac{1}{3}\,(P - x\,\partial_{x}),$ we
obtain
\begin{eqnarray*}
||(t\,\partial_{t})^{m}\,\partial_{x}^{l}u||_{H_{x}^{1}(I_{x_{0}})}
& = &
3^{-m}\,||(P - x\,\partial_{x})^{m}\,\partial_{x}^{l}u||_{H_{x}^{1}(I_{x_{0}})}   \\
& \leq & 3^{-m}\sum_{m=j_{1} + j_{2}}\frac{m!}{j_{1}!\,j_{2}!}\,
||(x\,\partial_{x})^{j_{1}}\,P^{j_{2}}\partial_{x}^{l}u||_{H_{x}^{1}(I_{x_{0}})} \\
& \leq & 3^{-m}\sum_{m=j_{1} + j_{2}}\frac{m!}{j_{1}!\,j_{2}!}\,
||(x\,\partial_{x})^{j_{1}}\,\partial_{x}^{l}(P - l)^{j_{2}}u||_{H_{x}^{1}(I_{x_{0}})} \\
& \leq & 3^{-m}\sum_{m=j_{1} + j_{2} +
j_{3}}\frac{m!}{j_{1}!\,j_{2}!\,j_{3}!}\,l^{j_{3}}\,
||(x\,\partial_{x})^{j_{1}}\,\partial_{x}^{l}P^{j_{2}}u||_{H_{x}^{1}(I_{x_{0}})}.
\end{eqnarray*}
where we replace $j_{2}$ into $j_{2} + j_{3}.$ Now, using the
induction hypothesis we have (with $B_{17}\geq A_{16}\,B_{16}$)
\begin{eqnarray}
\lefteqn{||(t\,\partial_{t})^{m}\,\partial_{x}^{l}u||_{H_{x}^{1}(I_{x_{0}})}} \nonumber \\
& \leq & 3^{-m}\,\sum_{m=j_{1} + j_{2} +
j_{3}}\frac{m!}{j_{1}!\,j_{2}!\,j_{3}!}\,l^{j_{3}}\,
c_{3}\,B_{17}^{j_{1} + j_{2} + l}\,(j_{1} + j_{2} + l)!\nonumber \\
\label{e467}& \leq & 3^{-m}\,c_{3}\,B_{17}^{m + l}\,(m +
l)!\sum_{m=j_{1} + j_{2} +
j_{3}}B_{17}^{-j_{3}}\frac{m!}{j_{1}!\,j_{2}!\,j_{3}!}\,l^{j_{3}}\,
\frac{(j_{1} + j_{2} + l)!}{(m + l)!},
\end{eqnarray}
Observing that
$\displaystyle
l^{j_{3}}\, \frac{(j_{1} + j_{2} + l)!}{(m + l)!}\leq 1,
$
we obtain in \eqref{e467}
\begin{eqnarray*}
||(t\,\partial_{t})^{m}\,\partial_{x}^{l}u||_{H_{x}^{1}(I_{x_{0}})}
& \leq  &
3^{-m}\,c_{3}\,(2 + B_{17}^{-1})^{m}\,B_{17}^{l + m}\,(l + m)!\\
& \leq & c_{4}\,A_{17}^{l + m}\,(l + m)!
\end{eqnarray*}
where we take $A_{17}=\mbox{max}\{B_{17},\,3^{-1}\,B_{17}\,(2 +
B_{17}^{-1})\}.$ We show that for some positive constants $c_{4},$
$A_{18}$ and $B_{18}$ we have
\begin{eqnarray}
\label{e468}||(t\,\partial_{t})^{j}\,\partial_{t}^{m}\partial_{x}^{l}u||_{H_{x}^{1}(I_{x_{0}})}\leq
c_{4}\,A_{18}^{j + m + l}\,B_{18}(j + m + l)!,\quad
j,\,l,\,m=0,\,1,\,2,\,\ldots
\end{eqnarray}
Induction in $m.$
\begin{eqnarray*}
||(t\,\partial_{t})^{j}\,\partial_{t}^{m +
1}\partial_{x}^{l}u||_{H_{x}^{1}(I_{x_{0}})} & \leq &
||\partial_{t}(t\,\partial_{t} - I)^{m}\,
\partial_{t}^{m}\partial_{x}^{l}u||_{H_{x}^{1}(I_{x_{0}})}\\
& = & t^{-1}\,||t\,\partial_{t}(t\,\partial_{t} - I)^{j}\,
\partial_{t}^{m}\partial_{x}^{l}u||_{H_{x}^{1}(I_{x_{0}})}\\
& \leq & (t_{0} - \epsilon)^{-1}\sum_{j_{1}=0}^{j}{j\choose
j_{1}}\;||(t\,\partial_{t})^{j_{1} + 1}\,
\partial_{t}^{m}\partial_{x}^{l}u||_{H_{x}^{1}(I_{x_{0}})}.
\end{eqnarray*}
Using the induction hypothesis
\begin{eqnarray*}
\lefteqn{||(t\,\partial_{t})^{j}\,\partial_{t}^{m + 1}\partial_{x}^{l}u||_{H_{x}^{1}(I_{x_{0}})}} \\
& \leq & (t_{0} - \epsilon)^{-1}\sum_{j_{1}=0}^{j}{j\choose j_{1}}\;
c_{4}\,A_{18}^{j_{1} + l + m + 1}\,B_{18}^{m}\,(j_{1} + l + m + 1)!\\
& = & c_{4}\,(t_{0} - \epsilon)^{-1}\,A_{18}^{j_{1} + l + m + 1}
\,B_{18}^{m}\,(j_{1} + l + m + 1)!\\
&  & \times\; \sum_{j_{1}=0}^{j}\frac{A_{18}^{-(j - j_{1}}}{(j -
j_{1})!}\,
{j\choose j_{1}}\;\frac{(j_{1} + m + l + 1)!\,(j - j_{1})!}{(j + m + l + 1)!}\\
& = & c_{4}\,(t_{0} - \epsilon)^{-1}\,e^{-A_{18}}\,A_{18}^{j + l + m
+ 1}
\,B_{18}^{m}\,(j + l + m + 1)!\\
& \leq & c_{4}\,\,A_{18}^{j + l + m + 1} \,B_{18}^{m}\,(j + l + m +
1)!
\end{eqnarray*}
where we take $B_{18}$ larger than $(t_{0} -
\epsilon)^{-1}\,e^{-A_{18}}.$ Finally, we choose $j=0$ in
\eqref{e468} and take $c_{2}=c_{4}$ and $A_{15}=A_{18}\,B_{18}.$ The
result of analyticity follows.

\section*{Acknowledgments} 

This work has been supported by
Fondap in Applied Mathematics (Project \# 15000001),
CNPq/CONICYT Project, \# 490987/2005-2 (Brazil) and
\# 2005-075 (Chile).

\end{document}